\documentstyle[12pt,amscd]{amsart}

\def\xnat{X^{\natural}}
\def\fm{Fourier-Mukai transform}
\def\p{{\mathcal{P}}}
\def\g{{\mathfrak{g}}}
\def\o{{\mathcal{O}}}
\def\f{{\mathcal{F}}}
\def\G{{\mathcal{G}}}
\def\d{{\mathcal{D}}}
\def\b{{\mathcal{A}}}
\def\for{\underline{for}}
\def\too{\longrightarrow}

\def\modsp{\text{Mod}(\o_X)_{sp}}
\def\modcxn{\text{Mod}(\o_Y)_{cxn}}
\def\obmodsp{\text{Ob Mod}(\o_X)_{sp}}
\def\obmodcxn{\text{Ob Mod}(\o_Y)_{cxn}}
\numberwithin{equation}{section}
\theoremstyle{plain}
\newtheorem{Def}{Definition}[section]
\newtheorem{Prop}[Def]{Proposition}
\newtheorem{Lem}[Def]{Lemma}

\newtheorem{Thm}[Def]{Theorem}
\newtheorem{Cor}[Def]{Corollary}
\newtheorem{Fund}{Mukai's Theorem}

\newtheorem{Funda}{Theorem Mukai 2}

\begin{document}

\title{Sheaves with connection on abelian varieties}
\author{Mitchell Rothstein}
\address{Department of Mathematics\\ University of Georgia\\ Athens, GA
30602}
\email{rothstei@@math.uga.edu}
\thanks{Research supported by NSF Grant No. 58-1353149}
\date{10/24/95}
\maketitle

\section{Introduction}\label{introduction}
\bigskip Let $X$ and $Y$ be abelian varieties over an algebraically closed
field
$k$, dual to one another, and let $\text{Mod}({\mathcal{O}}_X)$ and
$\text{Mod}({\mathcal{O}}_Y)$ be their respective categories of quasicoherent
$\o$-modules. Mukai proved  in \cite{Muk}  that the derived categories
$D\text{Mod}({\mathcal{O}}_X)$ and
$D\text{Mod}({\mathcal{O}}_Y)$ are equivalent  via a transform now known as the \fm,
\begin{equation}
{\mathcal{S}}_1({\mathcal{F}})=\alpha_{2*}(\p\otimes
\alpha_1^*(-1)^*({\mathcal{F}}))\ ,
\end{equation} where
$\p$ is the Poincar\'e sheaf and $\alpha_1$ and $\alpha_2$ are the
projections from $X\times Y$ to $X$ and $Y$ respectively.   A few years
earlier, Krichever
\cite{K} rediscovered a construction due originally  to Burchnall and
Chaundy \cite{BC}, by which the affine coordinate ring of a projective curve
minus a point  may be imbedded in the ring of formal differential operators
in one variable.  The construction involves the choice of a line bundle on
the curve, and Krichever took the crucial step of asking, in the case of a
smooth curve, how the imbedding varies when the line bundle moves linearly
on the Jacobian. The answer is now well-known, that the imbeddings satisfy
the system of differential equations known as the KP-hierarchy. In fact,
the Krichever construction is an instance of the
\fm, with the crucial addition that the transformed sheaf is not only an
${\mathcal{O}}_Y$-module but a ${\mathcal{D}}_Y$-module, where ${\mathcal{D}}_Y$ is the sheaf of
linear differential operators on
$Y$, \cite{N1} \cite{N2} \cite{R}.

This example serves as the inspiration for the present work, which is
concerned  with the role of the \fm\ in the theory of  sheaves on $Y$
equipped with a connection.  The main point is that in  the derived
category, all sheaves on $Y$ with connection are constructed by the \fm\ in
a manner directly generalizing the Krichever construction. The connection
need not be integrable, though the paper focuses mostly on that case.

The basic idea is the following.   Set
\begin{equation}
{\mathfrak{g}} = H^1(X, {\mathcal{O}})\ .
\end{equation} Then there is a tautological extension
\begin{equation}\label{taut seq} 0 \longrightarrow {\mathfrak{g}}^* \otimes {\mathcal{O}}
\longrightarrow
{\mathcal{E}}
\overset{\mu}{\longrightarrow} {\mathcal{O}} \longrightarrow 0
\end{equation} given by the extension class
$1 \in \text{End}({\mathfrak{g}}^*) = \text{Ext}^1({\mathcal{O}}, {\mathfrak{g}}^*
\otimes {\mathcal{O}})$.   Now let $\f$ be any quasicoherent sheaf of
$\o_X$-modules, and tensor the sequence \eqref{taut seq} with $\f$:
\begin{equation}\label{split it!} 0 \longrightarrow {\mathfrak{g}}^* \otimes \f
\longrightarrow
{\mathcal{E}} \otimes \f
\overset{\mu_{\f}}{\longrightarrow} {\mathcal{F}} \longrightarrow 0.
\end{equation} We refer to any splitting of sequence \eqref{split it!} as a
{\it splitting on $\f$.}  Let $\modsp$ denote the category
of pairs
$(\f,\psi)$, where $\f$ is a sheaf on $X$ and  $\psi:\f\too {\mathcal{E}}\otimes
\f$ is a splitting, with the obvious morphisms.

Let $\modcxn$ denote the category of quasicoherent sheaves on $Y$ equipped
with a connection.   In section \ref{mainthm} we use the \fm\ to establish
an equivalence of bounded derived  categories:
\begin{equation}\label{equiv1} D^b\modsp\ \ \leftrightarrow\ \ D^b\modcxn\ .
\end{equation} Note that   the extension class of
\eqref{split it!} belongs ${\mathfrak{g}}^*\otimes H^1(X,{\mathcal{E}}nd({\mathcal{F}}))$ and is
therefore a linear map
\begin{equation} H^1(X,\o)\longrightarrow H^1(X,{\mathcal{E}}nd({\mathcal{F}}))\ .
\end{equation} It is easily seen to be  the map on cohomology induced by the
map
${\mathcal{O}} \longrightarrow {\mathcal{E}}nd({\mathcal{F}})$.  Thus

\begin{Prop}\label{affine space} Given an $\o$-module $\f$, there is a
splitting on
${\mathcal{F}}$ if and only if the natural map
\begin{equation}
{\mathcal{O}} \longrightarrow {\mathcal{E}}nd({\mathcal{F}})
\end{equation} induces the $0$-map
\begin{equation}\label{map} H^1(X, {\mathcal{O}})\overset{0}{\longrightarrow} H^1(X,
{\mathcal{E}}nd({\mathcal{F}}))\ \ .
\end{equation} Moreover,  if (\ref{map}) holds, then the set of such
splittings is an affine space over \hbox{${\mathfrak{g}}^*\otimes H^0(X,\cal
End({\mathcal{F}}))$.}
\end{Prop} The intuitive idea behind the equivalence \eqref{equiv1} is the
following. Let
$g=\text{dim}(Y)$.  Let ${\mathcal{U}}_1,...,{\mathcal{U}}_k$ be an affine open cover of
$X$ and for $i=1,...,g$, let $\{c(i)_{m,n}\}\in Z^1(\{{\mathcal{U}}\},{\mathcal{O}})$ be a
1-cocycle, such that the classes $[c(1)],...,[c(g)]$ form a basis for
$H^1(X,{\mathcal{O}})$.   Let $\xi_1,...,\xi_g$ denote this basis, and let
$\omega^1,...,\omega^g$ denote the dual basis for
$H^1(X,{\mathcal{O}})^*$.  In light of proposition \ref{affine space}, a splitting
on $\f$ amounts to a collection of  endomorphisms \hbox{$\psi(i)_n\in
\Gamma({\mathcal{U}}_n,{\mathcal{E}}nd( (-1_X)^*({\mathcal{F}}))$} such that
\begin{equation}
\psi(i)_n-\psi(i)_m=\text{multiplication by}\ -c(i)_{n,m} \ \ .\label{cobdry}
\end{equation}   Let $\G={\mathcal{S}}_1(\f)$.  Then the collection of endomorphisms
$\psi=\{\psi(i)_n\}$ endows $\G$ with a connection in the following way. For
each
$n$, there is a connection
$\nabla_n$ relative to ${\mathcal{U}}_n$ on
${\mathcal{P}}\vert_{{\mathcal{U}}_n\times Y}$, such that on the overlaps,
\begin{equation}\label{overlaps}
\nabla_n-\nabla_m=c_{nm}\ \ .
\end{equation} Therefore, one gets a connection relative to $X$ on
$\alpha_1^*(-1_X)^*({\mathcal{F}})\otimes{\mathcal{P}}$ by defining
\begin{equation}\label{relcxn}
\nabla(\phi\otimes\sigma)=\phi\otimes\nabla_n(\sigma)+
\sum_i\omega^i\psi(i)_n(\phi)\otimes
\sigma
\end{equation} for $\phi\in(-1_X)^*({\mathcal{F}})$ and $\sigma\in{\mathcal{P}}$. Now one
applies
$\alpha_{2*}$ to produce  a connection on $\G$.

In the Krichever construction, $X$ is the  Jacobian  of a smooth curve $C$
with a base point $P$, and  ${\mathcal{F}}$ is
${\mathcal{O}}_C(*P)$, regarded as a sheaf on $X$ by the abel map.  The case where
$X$ is an arbitrary abelian variety and $\f$ is of the form
$\G\otimes\o_X(*D)$ for a coherent sheaf $\G$ and an ample hypersurface
$D\subset X$ has been studied in \cite{N1} and \cite{N2}.

Now consider the curvature tensor. To each object
$({\mathcal{F}},\psi)\in \obmodsp$ we associate a section
\begin{equation} [\psi, \psi]\in
\wedge^2 {\mathfrak{g}}^* \otimes  End({\mathcal{F}})\ ,
\end{equation} simply by taking the commutator
$[\psi(i)_n,\psi(j)_n]$, which, by \eqref{cobdry}, is independent of the
chart. Applying the Fourier-Mukai transform  to morphisms, one has
\begin{equation}
{\mathcal{S}}_1([{\psi,
\psi}])
\in
\wedge^2
{\mathfrak{g}}^* \otimes End({\mathcal{S}}_1(\f))\  .
\end{equation} Letting $[\nabla,\nabla]$ denote curvature, one has
(Proposition
\ref{curvature}),
\begin{equation}\label{curvature.intro}
{\mathcal{S}}_1([\psi,\psi])= [\nabla_{\psi},\nabla_{\psi}]\ .
\end{equation} In particular, ${\mathcal{S}}_1$ restricts to a functorial
correspondence
\begin{equation} (\f,\psi)\ \text{with}\ [\psi,\psi]=0\ \ \rightsquigarrow\
\
{\mathcal{D}}\text{-module structure on}\ {\mathcal{S}}_1(\f)\ .
\end{equation} We prove that this also induces an equivalence of bounded
derived categories (theorem \ref{equiv2}).

The main point  regarding the integrable case is the following. Let
\begin{equation} \xnat \overset{\pi}{\longrightarrow} X
\end{equation} denote the $\g^*$-principal bundle associated to the extension
${\mathcal{E}}$.  It is known that
$\xnat$ is the moduli space of line bundles on $Y$ equipped with an
integrable connection. For a discussion of $\xnat$ in greater generality, see
\cite{Me}, \cite{Ros}, and \cite{S}. Let
${\mathcal{A}}=\pi_*({\mathcal{O}}_{\xnat})$. Since $\pi$ is an affine morphism, the
category of ${\mathcal{O}}_{\xnat}$-modules is equivalent to category of $\cal
A$-modules. Then the subcategory of $\modsp$  whose objects satisfy $[\psi,
\psi]=0$ is precisely $\text{Mod}(\b)=\text{Mod}(\o_{\xnat})$ (Proposition
\ref{curvature}). Thus we have an equivalence of categories
\footnote{The referee informs us that this equivalence also appears in an
unpublished preprint by Laumon \cite{L}}
\begin{equation}\label{equiv2.intro} D^b\text{Mod}(\o_{\xnat})\ \
\leftrightarrow\
\ D^b\text{Mod}(\d_Y)\ .
\end{equation}

The outline of the paper is as follows.  In section \ref{mainthm} we prove
the basic equivalence theorem.   By way of illustration, section
\ref{Mat's thm} offers a new proof of a theorem of Matsushima on vector
bundles with a connection.  In section \ref{integrable case} the equality
\eqref{curvature.intro} is established
and
the  equivalence
\eqref{equiv2.intro} is proved.  Section \ref{some examples} gives some
examples.  In particular, we establish the formula
\begin{equation}
\hat{{\mathcal{A}}} = {\mathcal{D}}_{\{0\}\rightarrow Y}\ .
\end{equation}
Sections \ref{duality} and \ref{char var} contain  general results
about coherence, holonomicity and the characteristic variety.   Section
\ref{krich} illustrates the  theory in the case of the Krichever
construction.  In sections \ref{nak} and \ref{PDOs} we refine  and extend
several results of Nakayashiki on characteristic varieties of BA-modules
and commuting rings of  matrix partial differential operators.  This last
topic is a natural setting for the further study of integrable systems;
some brief remarks on this relationship are included at the end.  Further
applications to nonlinear partial differential equations will appear in a
future work.
\vskip 12pt

\noindent Acknowledgements:   The author wishes to thank Robert Varley for
many valuable discussions.

\section{First equivalence theorem}\label{mainthm}

We adopt the following sign conventions for the \fm:
\begin{align}
\text{Mod}(\o_X)&\overset{{\mathcal{S}}_1}\too\text{Mod}(\o_Y)\notag\\
{\mathcal{S}}_1({\mathcal{F}})&=\alpha_{2*}(\p\otimes
\alpha_1^*(-1)^*({\mathcal{F}}))\ ,
\end{align}
\begin{align}
\text{Mod}(\o_Y)&\overset{{\mathcal{S}}_2}\too\text{Mod}(\o_X)\notag\\
{\mathcal{S}}_2({\mathcal{G}})&=\alpha_{1*}(\p\otimes
\alpha_2^*({\mathcal{G}}))\ .
\end{align}
The fundamental result in \cite{Muk} is
\vskip.1in
\begin{Fund}  The derived functors
\begin{align} D^b\text{Mod}({\mathcal{O}}_X)&\overset{R{\mathcal{S}}_1}{\longrightarrow}
D^b\text{Mod}({\mathcal{O}}_Y)\notag\\  D^b\text{Mod}({\mathcal{O}}_Y) &\overset{R\cal
S_2}{\longrightarrow} D^b\text{Mod}({\mathcal{O}}_X)\end{align} are defined, and
\begin{align} R{\mathcal{S}}_1  R{\mathcal{S}}_2 &= T^{-g}\\ R{\mathcal{S}}_2
R{\mathcal{S}}_1 &= T^{-g}\ ,
\end{align} where $T$ is the shift autormophism on the derived category,
\begin{equation} (TF)^n = F^{n+1}.\notag
\end{equation}
\end{Fund}
\begin{pf}  \cite[p. 156]{Muk} \end{pf}

The other key result of [Muk] for our purposes is that the \fm\  exchanges
tensor product and Pontrjagin product.  Letting
$\beta_1$ and $\beta_2$ denote the projections from $Y \times Y$ to
$Y$ and denoting the group law on $Y$ by $m$,  the Pontrjagin product is
defined by
\begin{equation}
{\mathcal{F}} * {\mathcal{G}} = m_*(\beta^*_2({\mathcal{F}}) \otimes \beta^*_1({\mathcal{G}}))\ .
\end{equation}
\begin{Funda}
\begin{align}R{\mathcal{S}}_2({\mathcal{G}}_1
\overset{R}{\underset{=}{*}} {\mathcal{G}}_2) &= R{\mathcal{S}}_2({\mathcal{G}}_1)
\overset{L}{\underset{=}{\otimes}} RS_2({\mathcal{G}}_2)\label{exch2} \\
R{\mathcal{S}}_1 ({\mathcal{F}}_1
\overset{L}{\underset{=}{\otimes}} {\mathcal{F}}_2) &= T^g(R{\mathcal{S}}_1({\mathcal{F}}_1)
\overset{R}{\underset{=}{*}} R{\mathcal{S}}_1({\mathcal{F}}_2))\ .\label{exch1}
\end{align}
\end{Funda}
\begin{pf}  \cite[p. 160]{Muk} \end{pf}

We want to apply this to the following situation.
Let ${\mathcal{I}} \subset \o_Y$ be the ideal
sheaf of $0 \in Y$, and let $k(0)$ be the skyscraper sheaf at the origin
with fiber $k$.  Then
$R^0{\mathcal{S}}_2(k(0))=\o_X$ and $R^i{\mathcal{S}}_2(k(0))=0$ for $i>0$.  Thus the
\fm\ takes the short exact sequence
\begin{equation}\label{nptext} 0 \longrightarrow \g^* \otimes_k k(0)
\longrightarrow
\o/{\mathcal{I}}^2
\longrightarrow k(0) \longrightarrow 0
\end{equation}
to a short exact sequence of vector bundles on $X$, and it is easy to see
that it is precisely the sequence \eqref{taut seq}.
Thus by theorem Mukai 2, if ${\mathcal{F}}
\overset{\psi}{\longrightarrow} {\mathcal{E}} \otimes {\mathcal{F}}$ is a splitting on
${\mathcal{F}}$, $\psi$ induces a morphism
\begin{equation} R{\mathcal{S}}_1({\mathcal{F}}) \overset{{\mathcal{S}}_1(\psi)}{\longrightarrow}
T^{-g}({\mathcal{O}}/{\mathcal{I}}^2) * R{\mathcal{S}}_1({\mathcal{F}})\ .
\end{equation} (We identify objects in an abelian category with
complexes concentrated in degree 0.)  If we now take $g^{th}$ cohomology, we
get
\begin{equation}\label{The point}
{\mathcal{S}}_1({\mathcal{F}}) \overset{H^g{\mathcal{S}}_1(\psi)}{\longrightarrow}
{\mathcal{O}}/{\mathcal{I}}^2 *
{\mathcal{S}}_1({\mathcal{F}}).
\end{equation} The point is this.  If ${\mathcal{G}}$ is any
$\o_Y$-module, there is a prolongation sequence
\begin{equation}\label{c} 0 \longrightarrow \Omega^1 \otimes_{{\mathcal{O}}}{\mathcal{G}}
\longrightarrow j({\mathcal{G}})
\overset{\nu_{\G}}{\longrightarrow} {\mathcal{G}} \longrightarrow 0,
\end{equation} such that a splitting of $\nu_{{\mathcal{G}}}$ is precisely a
connection on
${\mathcal{G}}$.
As a sheaf of abelian groups,
\begin{equation} j({\mathcal{G}}) = {\mathcal{G}} \oplus (\Omega^1 \otimes_{\o} {\mathcal{G}})\ ,
\end{equation} with $\o$-module structure
\begin{equation} f  (\phi, \omega \otimes \psi) =(f\phi, f \omega \otimes
\psi + df \otimes \phi)\ \ .
\end{equation} Thus a connection on ${\mathcal{G}}$ is a splitting of \eqref{c}.
Since $Y$ is an abelian variety, there is a characterization of
$j({\mathcal{G}})$ in terms of the Pontrjagin product.

\begin{Lem} \label{pont} For any $\o_Y$-module ${\mathcal{G}}$,
\begin{equation} j({\mathcal{G}}) = (\o/{\mathcal{I}}^2) * {\mathcal{G}}\ .
\end{equation}
\end{Lem}
\begin{pf} Let $Y_1\subset Y\times Y$
denote the first order neighborhood of
the diagonal, and let
$\pi_i: Y_1
\longrightarrow Y$, $i = 1, 2$, denote the two projections. Then
\begin{equation}
j(\f)=\pi_{2*}\pi_1^*(\f)\ .
\end{equation}
(This holds for any variety.)
Let $\tilde{Y} = \text{Spec}(\o/{\mathcal{I}}^2)$, the first order neighborhood of
$0$ in $Y$. Then $Y_1$ may be identified with
$Y \times \tilde{Y}$  in such a way that  $\pi_1$ corresponds to
projection onto the first factor and $\pi_2$ corresponds to the group law,
$\tilde m$. Let $\iota:\tilde Y\to Y$ denote the inclusion map.  Then
\begin{align}
\o/{\mathcal{I}}^2&*{\mathcal{G}}=
m_*(\beta_2^*\iota_*(\o_{\tilde Y})\otimes\beta_1^*({\mathcal{G}}))\notag\\
&=
m_*(1\times\iota)_*(\tilde\beta_2^*(\o_{\tilde Y})\otimes
\tilde\beta_1^*({\mathcal{G}}))\notag\\
&=\tilde m_*
\tilde\beta_1^*({\mathcal{G}}))=j({\mathcal{G}})\ .
\end{align}
\end{pf}
Combining this lemma with the map \eqref{The point}, we see that a
splitting on
${\mathcal{F}}$ induces a splitting of the prolongation sequence of ${\mathcal{S}}_1(\cal
F)$, i.e., a connection on ${\mathcal{S}}_1({\mathcal{F}})$.  So we have a functor
\begin{equation}
\modsp \overset{S_1}{\longrightarrow}
\modcxn\ .
\end{equation} We will check later that this description is equivalent to
the one given in the introduction.

Conversely, if we apply ${\mathcal{S}}_2$ to a splitting of the prolongation
sequence,
${\mathcal{G}}\overset{\tau}{\longrightarrow} j({\mathcal{G}}) =
{\mathcal{O}}/{\mathcal{I}}^2 * {\mathcal{G}}$,  theorem Mukai 2 gives a splitting
\begin{equation}
{\mathcal{S}}_2({\mathcal{G}}) \overset{\psi}{\longrightarrow} {\mathcal{E}} \otimes
{\mathcal{S}}_2({\mathcal{G}})\ .
\end{equation} So we have
\begin{equation}
\modcxn \overset{S_2}{\longrightarrow}
\modsp\ .
\end{equation}

The categories $\text{Mod}(\o_X)_{sp}$ and
$\text{Mod}(\o_Y)_{cxn}$ are abelian.  Moreover, objects in either
$\text{Mod}(\o_X)_{sp}$ or $\text{Mod}(\o_Y)_{cxn}$ may be resolved by a \v
Cech resolution with respect to an affine open cover of $X$ or $Y$.  Thus
the derived functors
\begin{equation} D^b\text{Mod}(\o_X)_{sp} \overset{RS_1}{\longrightarrow}
D^b\text{Mod}(\o_Y)_{cxn}
\end{equation}
\begin{equation} D^b\text{Mod}(\o_Y)_{cxn} \overset{RS_2}{\longrightarrow}
D^b\text{Mod}(\o_X)_{sp}
\end{equation} exist.

The main result of this section is

\begin{Thm}\label{equiv}
\begin{align} RS_1  RS_2 &= T^{-g} , \\ RS_2  RS_1 &= T^{-g} .
\end{align}
\end{Thm}

\begin{pf}  Let  $\zeta$ denote the functor $T^g  RS_1  RS_2$. Let
$\for$ denote the forgetful functor from
$D^b\text{Mod}({\mathcal{O}}_Y)_{cxn}$ to $D^b\text{Mod}({\mathcal{O}}_Y)$.
Then
\begin{equation}\for\ \zeta = \for
\end{equation}
by
Mukai's theorem.  In particular, for any object
$({\mathcal{F}},\nabla)
\in \obmodcxn$, $H^i\zeta({\mathcal{F}},\nabla) = 0$  for
$i > 0$. Thus
$\zeta(\f,\nabla)=(\f,\nabla')$ for some new connection $\nabla'$.

Let ${\mathcal{F}} \overset{\tau}{\longrightarrow} j({\mathcal{F}})$ denote the splitting
associated to $\nabla$, and let ${\mathcal{F}}
\overset{\tau'}{\longrightarrow} j({\mathcal{F}})$ denote the splitting associated
to $\nabla'$.  Let $ \psi$ denote the corresponding splitting on
${\mathcal{S}}_2(\f)$.  Then
$\psi$ is the
$0^{th}$ cohomology of
\begin{equation} R{\mathcal{S}}_2({\mathcal{F}}) \overset{R\cal
S_2(\tau)}{\longrightarrow} {\mathcal{E}}
\overset{L}{\underset{=}{\otimes}} R{\mathcal{S}}_2({\mathcal{F}})\ ,
\end{equation} from which it follows that $\tau'$ is the $0^{th}$ cohomology
of
\begin{equation}
\begin{CD}
{\mathcal{F}} @>{T^gR{\mathcal{S}}_1 R{\mathcal{S}}_2(\tau)}>> j({\mathcal{F}})\ .
\end{CD}
\end{equation} Thus $\tau = \tau'$, again by Mukai's Theorem.

Similarly, if $({\mathcal{G}}, \psi)$ is an object in $\text{Mod}({\mathcal{O}}_X)_{sp}$,
\begin{equation} T^gRS_2 RS_1({\mathcal{G}}, \psi) = ({\mathcal{G}}, \psi)\ .
\end{equation}
The next lemma then completes the proof
of the theorem.
\end{pf}

\begin{Lem} \label{bootstrap} Let  $C_1$ and $C_2$  be abelian categories,
and let
\begin{equation} D^bC_1 \underset{F_2}{\overset{F_1}{\rightarrow}}
D^bC_2
\end{equation}
be $\delta$-functors.  If $F_1$ and $F_2$ are isomorphic when restricted
to the subcategory $C_1\subset D^bC_1$, then they are isomorphic.
\end{Lem}
\begin{pf} This follows by induction on the cohomological length of an object
in the bounded derived category, using \cite[lemme 12.6, p.104]{Bo}  and the
triangle axiom TR3, \cite[p.28]{Bo}, once it is noted that the
constructions used there are functorial.
\end{pf}
\noindent{\bf Remark}\ \ Let
$\gamma_i$, $\gamma_{i,j}$ denote the projections on $X\times Y\times Y$.
The key to Mukai's theorem is the elementary formula
\begin{equation} \label{key}
\gamma^*_{1,2} ({\mathcal{P}}) \otimes \gamma^*_{1,3}({\mathcal{P}}) = (1\times
m)^*({\mathcal{P}})\ .
\end{equation}
This formula also plays a crucial but hidden role in theorem \ref{equiv},
which we would like to make explicit.

Let $\tilde{Y} = \text{Spec}(\o/{\mathcal{I}}^2)$.  Then $\p|_{X \times
\tilde{Y}}$ is a line bundle on $X \times \tilde{Y}$, trivial on $X \times
\{0\}$.  Set
$\tilde{\p} = \p|_{X \times \tilde{Y}}$, and let
$\tilde{\alpha}_1: X \times \tilde{Y} \rightarrow X$ be the projection.
Then ${\mathcal{E}} = \tilde{\alpha}_{1*}(\tilde{\p})$ and $\mu: {\mathcal{E}} \rightarrow
\o$ is the morphism which restricts a section of $\tilde{\p}$ to $X \times
\{0\}$.
Let
$\tilde\gamma_i$, $\tilde\gamma_{i,j}$ denote the projections on $X\times
Y\times \tilde Y$.  Then we get an infinitesimal form of \ref{key},
\begin{equation}\label{key2}
\tilde{\gamma}^*_{1,2} ({\mathcal{P}}) \otimes
\tilde{\gamma}^*_{1,3}(\tilde{{\mathcal{P}}}) = (1\times \tilde m)^*({\mathcal{P}}) \ .
\end{equation}
If ${\mathcal{G}}$ is a sheaf on $X\times Y$, then a connection on
${\mathcal{G}}$ relative to $X$ is an
isomorphism
\begin{equation}
\tilde{\gamma}_{1,2}^*({\mathcal{G}}) \approx (1\times \tilde m)^*({\mathcal{G}})
\end{equation}  restricting to the identity on $X \times Y$. Thus
\eqref{key2} says that
$\tilde{\gamma}^*_{1,3}(\tilde{{\mathcal{P}}})$ is the obstruction to endowing
${\mathcal{P}}$ with a connection relative to $X$.   Given
a sheaf $\f$ of $\o_X$-modules, a splitting on $\f$ is
precisely what is needed to cancel this
obstruction.  Indeed, a splitting may be regarded as an isomorphism
\begin{equation}
\tilde{\alpha}_1^*({\mathcal{F}})\overset{\psi}{\to}
\tilde{\p}\otimes\tilde{\alpha}_1^*({\mathcal{F}})
\end{equation}
restricting to the identity on $X$.
Applying $(-1_X,-1_{\tilde Y})$, we get
\begin{equation}
\tilde\alpha_1^*(-1)^*({\mathcal{F}})\longrightarrow
\tilde{\p}\otimes\tilde\alpha_1^*(-1)^*({\mathcal{F}})\ \ .
\end{equation}
If we then apply $\tilde{\gamma}_{1,2}^*$ to
$\p\otimes\alpha_1^*(-1)^*(\f)$ we find
\begin{align}
\tilde\gamma_{1,2}^*(\p\otimes&\alpha_1^*(-1)^*({\mathcal{F}}))=
\tilde\gamma_{1,2}^*(\p)\otimes\tilde\gamma_{1,3}^*
\tilde\alpha_1^*(-1)^*({\mathcal{F}})\notag\\  &\approx
\tilde\gamma_{1,2}^*(\p)\otimes\tilde\gamma_{1,3}^*(\tilde{\p})
\otimes\tilde\gamma_{1,3}^*\tilde\alpha_1^*(-1)^*({\mathcal{F}})
\notag\\
&=(1\times\tilde
m)^*(\p)\otimes\tilde\gamma_{1,3}^*\tilde\alpha_1^*(-1)^*({\mathcal{F}})\notag\\
&=(1\times\tilde m)^*(\p\otimes\alpha_1^*(-1)^*({\mathcal{F}}))\ \ ,
\label{cxnalt}
\end{align}
which is a relative connection on $\p\otimes\alpha_1^*(-1)^*(\f)$.   Then
apply $\alpha_{2*}$ to get a sheaf with connection on $Y$, and this is our
functor $S_1$.

\section{Matsushima's Theorem}\label{Mat's thm}

As an application of theorem
\ref{equiv}, we will give a new proof of Matsushima's Theorem
on the homogeneity  of vector bundles admitting a connection \cite{Mat}.
The key is the following lemma, which is of interest in
its own right.  Here we assume $char(k)=0$.

\begin{Lem} \label{fin supp}  Let $\f$ be a coherent $\o_X$-module with a
splitting.  Then $\f$ is finitely supported.
\end{Lem}

\begin{pf}  The proof is similar to that of lemma 3.3 in \cite{Muk}.
Assuming $\text{dim(supp}(\f))>0$, let $C$ be a curve contained in
$\text{supp}(\f)$, and let $\tilde C\overset{\pi}\too C$ be its
normalization.  Let $\f'=\pi^*(\f)$.  Then we get a non-zero vector bundle
$\f''$ on
$\tilde C$ upon taking the quotient of  $\f'$ by its torsion part.  Let
${\mathcal{E}}'$ denote the pullback of ${\mathcal{E}}|_C$ to $\tilde C$.  Then any
splitting on $\f$ induces a splitting
\begin{equation}
\f'\overset{\psi'}\too {\mathcal{E}}'\otimes \f'\ .
\end{equation}
Since $Tor( {\mathcal{E}}'\otimes \f')=
{\mathcal{E}}'\otimes Tor(\f')$, we get a splitting of the sequence
\begin{equation}\label{can't split}
0 \longrightarrow {\mathfrak{g}}^* \otimes \f''
\longrightarrow
{\mathcal{E}}' \otimes \f''
\longrightarrow {\mathcal{F}}'' \longrightarrow 0.
\end{equation}
But the extension class of the sequence \eqref{can't split} is the bottom
arrow in the commutative diagram
\begin{equation}
\begin{matrix} H^1(\tilde{C}, {\mathcal{O}}) \\ a\nearrow \qquad\qquad \searrow b
\quad \\ H^1(X, {\mathcal{O}}) \overset{e}{\longrightarrow} H^1(\tilde{C},
{\mathcal{E}}nd({\mathcal{F}}''))
\end{matrix}\ ,
\end{equation}
where $b$ is induced by the natural inclusion $\o\overset{\beta}\too \cal
End(\f'')$ and
$a$ is the derivative of the natural map $Pic(X)\to Pic(\tilde C)$.  In
particular, $a$ is not the $0$-map.  Moreover, in characteristic $0$,
$\beta$ splits by the trace, so $b$ is injective.  This is a contradiction,
for now $e\ne 0$ so that \eqref{can't split} does not split. \end{pf}

We then have
\begin{Thm}[Matsushima] Any vector bundle on $Y$ admitting a connection is
homogeneous.
\end{Thm}

\begin{pf}  Let $\G$ be such a vector bundle.  Then the cohomology sheaves
$R^i{\mathcal{S}}_2(\G)$ are $\o_X$-coherent and admit splittings.  By lemma
\ref{fin supp}, they are all finitely supported.  As in
\cite[example 3.2, p. 158]{Muk}, we then have $R^i{\mathcal{S}}_2(\G)=0$ for
$i\ne g$, and $\G$ is then the \fm\ of a finitely supported sheaf.  By
\cite[3.1]{Muk}, $\G$ is homogeneous. \end{pf}

\noindent{\bf Remark}\ \ The converse of the statement above is also part
of Matsushima's theorem.  That the converse can be proved by the \fm\
is already noted in \cite[prop 5.9]{N2}.

\section{Curvature tensor and the integrable case}\label{integrable case}

If $({\mathcal{G}},\nabla)$ is a sheaf with connection on $Y$, then its curvature
tensor is a linear map
\begin{equation}
[\nabla,\nabla]:\wedge^2(\g)\longrightarrow End({\mathcal{G}})\ \ .
\end{equation}
Before explaining how the curvature can be read off from the
transform of $\G$, we want to show that the functor $S_1$ has the
\v Cech description given in the introduction.
Let $\psi_1$ and $\psi_2$ be two splittings on a given
sheaf of $\o_X$-modules $\f$, and let $\eta=
\psi_1-\psi_2$.  By proposition \ref{affine space}, $\eta$ is a map
\begin{equation}
\eta:\g\too End(\f)\ .
\end{equation}
Applying the \fm\  to $End(\f)$, we get
\begin{equation}
{\mathcal{S}}_1(\eta):\g\too End({\mathcal{S}}_1(\f))\ .
\end{equation}
Denoting the two connections on ${\mathcal{S}}_1(\f)$ by $\nabla_1$  and
$\nabla_2$, it is easy to check
that
\begin{equation}\label{compat}
\nabla_2=
\nabla_1-{\mathcal{S}}_1(\eta)\ .
\end{equation}
Now let ${\mathcal{U}}\overset{\iota}\too X$ be an affine open subset,
equipped with a section $\rho\in\Gamma({\mathcal{U}},{\mathcal{E}})$ such that
$\mu(\rho)=1$.  Then $(\iota_*(\o_{{\mathcal{U}}}),\rho)$ is an object in
$\modsp$.  The corresponding object in $\modcxn$ is  a
connection on $\alpha_{2*}(\p|_{{\mathcal{U}}\times Y})$.  Since functions on
${\mathcal{U}}$ act as endomorphisms of $(\iota_*(\o_{{\mathcal{U}}}),\rho)$, this
connection is linear over $\alpha_{2*}\alpha_1^{-1}(\o_{{\mathcal{U}}})$.  So we in
fact have a connection relative to ${\mathcal{U}}$ on $\p|_{{\mathcal{U}}\times Y}$.  Call
this connection $\nabla^{\rho}$.   If $\f$ is any $\o_X$-module, $\rho$
induces a splitting on $\f|_{{\mathcal{U}}}$.  Then if $\psi$ is any other
splitting on $\f$, it must take the form
\begin{equation}
\psi(\kern .5em\cdot\kern .5em)=\kern .5em
\cdot
\kern .5em\otimes\rho+
\sum\limits_{i=1}^g\omega^i
\psi(i)^{\rho}(\kern .5em\cdot\kern .5em)\ ,
\end{equation}
for some endomorphisms $\psi(i)^{\rho}$.  Again the corresponding sheaf
with connection on $Y$ is the direct image of a sheaf with relative
connection on $X\times Y$, and by \eqref{compat} it has the
form of \eqref{relcxn}.

It is now easy to read off the curvature tensor of
$S_1({\mathcal{F}},\psi)$.
There is an exact sequence
\begin{align}
 0 \longrightarrow \wedge^2{\mathfrak{g}}^* \otimes {\mathcal{O}} &\longrightarrow
\wedge^2 {\mathcal{E}} \overset{\mu_2}{\longrightarrow} {\mathfrak{g}}^*
\otimes {\mathcal{O}} \longrightarrow 0,\notag\\
\mu_2(\rho_1 \wedge \rho_2) &= \mu(\rho_1) \rho_2 - \mu(\rho_2) \rho_1.
\end{align}
If we iterate $\psi$ and then skew-symmetrize, we get a map
\begin{equation}
{\mathcal{F}} \overset{[\psi, \psi]}{\longrightarrow} \wedge^2({\mathcal{E}})
\otimes{\mathcal{F}}.
\end{equation} One sees easily that $(\mu_2 \otimes 1) \circ [\psi,
\psi] = 0$, so
\begin{equation}
[\psi,
\psi]
\in   \wedge^2 {\mathfrak{g}}^*
\otimes End({\mathcal{F}})\ .\end{equation}
As noted in the introduction,
in terms of the family of endomorphisms $\psi(i)_n$, one
simply has
\begin{equation}\label{bracket}
[\psi, \psi]=\sum \omega^i\wedge\omega^j[\psi(i)_n,\psi(j)_n]\ .
\end{equation}
Applying the functor ${\mathcal{S}}_1$ to morphisms, we get
${\mathcal{S}}_1([\psi,
\psi])
\in
\wedge^2
{\mathfrak{g}}^* \otimes End({\mathcal{S}}_1(\f))$.
Set $S_1({\mathcal{F}},\psi)=({\mathcal{S}}_1({\mathcal{F}}),\nabla_{\psi})$.  Let
$[\nabla,\nabla]$ denote curvature.
\begin{Prop}\label{curvature}  $[\nabla_{\psi},\nabla_{\psi}]=\cal
S_1([\psi,\psi])$.
\end{Prop}

\begin{pf} First we claim that the relative connections $\nabla_n$ on
$\p\vert_{{\mathcal{U}}_n\times Y}$ are integrable.  Indeed, from \eqref{overlaps},
the curvature tensor of $\nabla_n$ is independent of the index $n$.
To show
that it is 0, we have only to consider the case where $0 \in {\mathcal{U}}_n$
and note
that all connections on the trivial bundle ${\mathcal{O}}_Y$ are integrable.

Denoting  the relative connection on $\p\otimes\alpha^*_1(-1_X)^*({\mathcal{F}})$
by
$\nabla=\sum \omega^i \nabla(i)$,
\begin{align}
[\nabla(i), \nabla(j)]&=[\nabla(i)_n\otimes 1 - 1\otimes\psi(i)_n,
\nabla(j)_n\otimes 1 - 1\otimes\psi(j)_n]\notag\\
&=1\otimes[\psi(i)_n,
\psi(j)_n]\ \ .
\end{align}
\end{pf}

Similarly, given $({\mathcal{G}},\nabla)$ in $\obmodcxn$,  we have $\cal
S_2([\nabla,\nabla])\in
\wedge^2({\mathfrak{g}}^*)\otimes End({\mathcal{S}}_2({\mathcal{G}}))$.  Setting
$S_2({\mathcal{G}},\nabla)=({\mathcal{S}}_2({\mathcal{G}}),\psi_{\nabla})$,
one proves
\begin{Prop}  $[\psi_{\nabla},\psi_{\nabla}]={\mathcal{S}}_2([\nabla,\nabla])$.
\end{Prop}

Turning now to the integrable case,
consider the group extension
\begin{equation}
0 \longrightarrow {\mathfrak{g}}^* \longrightarrow
X^{\natural}
\overset\pi{\longrightarrow} X \longrightarrow 0
\end{equation}
mentioned in the introduction.
The morphism $\pi$ being affine, $\xnat$ is characterized by $\pi_*$ of its
structure sheaf.  Define a sheaf
${\mathcal{A}}$ of $\o_X$-modules as follows.
For each affine open ${\mathcal{U}} \subset X$ and each
$\rho
\in
{\mathcal{E}}({\mathcal{U}})$ such that $\mu(\rho) = 1$, introduce independent variables
$x^{\rho}_1, \cdots, x^{\rho}_g$.  Then
\begin{equation}
{\mathcal{A}}|_{{\mathcal{U}}}=
\o|_{{\mathcal{U}}}[x^{\rho}_1, \cdots, x^{\rho}_g]\ .
\end{equation}
Glue these sheaves together by the rule that if
$\tilde{\rho} =\rho + \operatornamewithlimits\sum\limits^g_{i=1}
\omega^i f_i$, then
\begin{equation}
 x^{\rho}_i = x^{\tilde{\rho}}_i + f_i
\end{equation} as sections of ${\mathcal{A}}$.
\begin{Def}   $\xnat=\text{Spec}({\mathcal{A}})$.
\end{Def}

\begin{Prop}\label{subcat} The full subcategory of
$\modsp$ whose objects $({\mathcal{F}}, \psi)$ satisfy $[\psi,
\psi] = 0$ is canonically isomorphic to the category of $\o$-modules
on $\xnat$, which is to say the category of ${\mathcal{A}}$-modules on $X$.
\end{Prop}

\begin{pf}  Let
$\rho$ be as above.  Since ${\mathcal{A}}$ is locally  ${\mathcal{O}}[x^{\rho}_1,
\dots, x^{\rho}_g]$, to give an ${\mathcal{A}}$-module structure to an $\cal
O_X$-module ${\mathcal{F}}$ is to choose, for every $\rho$, a set of commuting
endomorphisms $\psi^{\rho}_1,
\dots, \psi^{\rho}_g \in\Gamma({\mathcal{U}},{\mathcal{E}}nd({\mathcal{F}}))$, such that if
\begin{equation}
\tilde{\rho} = \rho + \sum \omega^i f_i\ ,
\end{equation}
then $\psi^{\rho}_i - \psi^{\tilde{\rho}}_i=$
multiplication by
$f_i$.  This is clearly the same as lifting ${\mathcal{F}}$ to an object
$({\mathcal{F}}, \psi) \in \obmodsp$ such that $[\psi, \psi] = 0$.
\end{pf}

Now we may restrict our functors $S_1$ and $S_2$ to the subcategories
$\text{Mod}(\o_{\xnat})\subset\modsp$ and
$\text{Mod}({\mathcal{D}}_Y)\subset\modcxn$ respectively, to get  functors
\begin{align}
\text{Mod}(\o_{\xnat})&\overset{S_1}\longrightarrow \text{Mod}(\cal
D_Y)\notag\\
\text{Mod}({\mathcal{D}}_Y)&\overset{S_2}\longrightarrow\text{Mod}(\o_{\xnat})\ \ .
\end{align}

We then get our second equivalence theorem, whose proof is the same as that
of the first.

\begin{Thm}\label{equiv2}
The derived functors
\begin{equation} D^b\text{Mod}(\o_{\xnat})
\overset{RS_1}{\longrightarrow} D^b\text{Mod}({\mathcal{D}}_Y)
\end{equation}
\begin{equation} D^b\text{Mod}({\mathcal{D}}_Y)
\overset{RS_2}{\longrightarrow} D^b\text{Mod}(\o_{\xnat})
\end{equation} exist, and satisfy
\begin{align} RS_1  RS_2 &= T^{-g} , \notag\\ RS_2  RS_1 &= T^{-g} .\notag
\end{align}
\end{Thm}

\section{Some examples}\label{some examples}

\subsection{}  Using theorem \ref{equiv2}, one easily recovers the
observation that $\xnat$ is the moduli space of degree-0 line bundles on $Y$
equipped with a connection \cite{Me}.
Indeed, let $\sigma:Z\to \xnat$ be any morphism and let
$\sigma'=\pi\circ\sigma$.   Let $\p_{\sigma'}$ be the  line bundle
on $Z\times Y$ induced by $\sigma'$.  Then the ${\mathcal{D}}_Y$-module
$S_1(\sigma_*(\o_Z))$
is in fact the direct image of a flat connection relative to $Z$ on
$\p_{\sigma'}$.
In particular, we have the 1:1 correspondence
$$\text{ points of}\ \xnat \leftrightarrow \text{line bundles with
connection on}\ Y\ .$$

\subsection{} Let $J\subset\text{Sym}(\g)$ be
an ideal.  Then $\text{Spec}(\text{Sym}(\g)/J)$ is a subvariety of the
fiber $\pi^{-1}(0)$, and we may apply $S_1$ to its structure sheaf.
This gives a ${\mathcal{D}}$-module
structure on $\o_Y\otimes_k\text{Sym}(\g)/J$.  Writing sections of ${\mathcal{D}}$
in the standard format with functions on the left and invariant
differential operators on the right, make the identification
\begin{equation}
{\mathcal{D}}\simeq\o_Y\otimes_k\text{Sym}(\g)\ .
\end{equation}
Then $\o_Y\otimes_kJ$ is a left ideal, and $S_1(
k(0)\otimes_k\text{Sym}(\g)/J)$ is the quotient ${\mathcal{D}}$-module.  In
particular,
\vskip 12pt
\noindent ${\mathcal{D}}_Y$ is the \fm\ of the structure sheaf of the subgroup
${\mathfrak{g}}^* =\pi^{-1}(0)\subset\xnat$.

\subsection{}\label{B itself}   Going the other way, we may
ask which $\d_Y$-module transforms to ${\mathcal{A}}$.  By theorem \ref{equiv2}
this is the same as asking for the \fm\ of
${\mathcal{A}}$.  Observe first that there is an important filtration $\{\cal
A(m)\}$ on
${\mathcal{A}}$ coming from  the local identifications of
${\mathcal{A}}$ with a polymonial algebra over ${\mathcal{O}}$.   One has
exact sequences
\begin{equation}\label{filtration}
 0 \rightarrow {\mathcal{A}}(m) \rightarrow {\mathcal{A}}(m+1)
\rightarrow Sym^{m+1}({\mathfrak{g}}) \otimes_k {\mathcal{O}}_X \rightarrow 0\ \ .
\end{equation}
We then have
\begin{Lem}
$R^i{\mathcal{S}}_1({\mathcal{A}})=0$ for $i\ne g$, and
$R^g{\mathcal{S}}_1({\mathcal{A}})$ is supported at the origin.\end{Lem}\begin{pf}
By \eqref{filtration} the question reduces to the corresponding assertion
about
$\o$.  The latter is proved in \cite[p.126]{Mum}.\end{pf}

It remains to identify the ${\mathcal{D}}_Y$-module $R^gS_1({\mathcal{A}})$.
If $\G$ is an $\o_Y$-module, we may form $\d
\otimes_{\o} \G$, using right multiplication by functions as the
$\o$-module structure on $\d$.  Then $\d \otimes_{\o} \G$ has a left
$\d$-module structure.  In particular,  define
\begin{equation}
{\mathcal{D}}_{\{0\}\rightarrow Y}\overset{def}=\d\otimes k(0)\ .
\end{equation}
\begin{Prop}\label{fm of B}  $R^gS_1({\mathcal{A}}) = {\mathcal{D}}_{\{0\}\rightarrow Y}$.
\end{Prop}
\begin{pf*}{First proof}  By theorem \ref{equiv2}, it suffices to prove that
$R^0S_2({\mathcal{D}}_{\{0\}\rightarrow Y}) = {\mathcal{A}}$.  Explicitly, $\cal
D_{\{0\}\rightarrow Y}$ is the sheaf of $k$-vector spaces
$\text{Sym}({\mathfrak{g}})(0)$, with the following ${\mathcal{D}}$-module structure.
Given $L \in \text{Sym}({\mathfrak{g}})$ and $\xi \in {\mathfrak{g}}$,
$\xi$ acts on $L$ by $\xi \cdot L = \xi L$.   To give the action of
${\mathcal{O}}_Y$, we need only consider $f \in {\mathcal{O}}_{Y,0}$.  Expand the
differential operator
$f L$   in the form
$\sum \xi^I g_I$,\ $g_I \in {\mathcal{O}}_{Y, 0}$.  Define
$$ :f L: = \sum  \xi^I g_I(0)\in\text{Sym}({\mathfrak{g}})\ .
$$ Then $f \cdot L = :f L:$\ .

Let ${\mathcal{F}} = R^0S_2({\mathcal{D}}_{\{0\}\rightarrow Y})$.  We want to prove that
${\mathcal{F}}$ is a free rank-1 ${\mathcal{A}}$-module.  We have the global section
$1 \in\text{Sym}({\mathfrak{g}}) = H^0(Y, {\mathcal{D}}_{\{0\} \rightarrow
Y})$.
Let ${\mathcal{U}} \subset X$ be an open subset and let
$\sigma$ be a meromorphic section of ${\mathcal{P}}$ on ${\mathcal{U}} \times Y$ such that
$\sigma(0) = 1$, where $\sigma(0) = \sigma|_{{\mathcal{U}}
\times \{0\}}$.  Since $\sigma$ is regular on the support of
$\alpha^*_2({\mathcal{D}}_{\{0\}\rightarrow Y})$, there is a well-defined section
$$
\sigma\otimes 1 \in \Gamma({\mathcal{U}} \times Y, {\mathcal{P}} \otimes
\alpha^*_2({\mathcal{D}}_{\{0\}\rightarrow Y})) = \Gamma({\mathcal{U}}, {\mathcal{F}})\ .
$$ If $f$ is a meromorphic function on ${\mathcal{U}} \times Y$, regular on
${\mathcal{U}}
\times \{0\}$, then
$$ f \sigma \otimes 1 = \sigma \otimes :f 1  :\ = f(0) \sigma \otimes 1\ .
$$ Thus $\sigma \otimes 1$ is independent of $\sigma$, and defines a
global section of $\f$.  Now let
$\rho
\in \Gamma({\mathcal{U}},
{\mathcal{E}})$ satisfy $\mu(\rho) = 1$.  Then $\rho$ determines both a relative
connection $\nabla$ on ${\mathcal{P}}|_{{\mathcal{U}} \times Y}$ and a set of sections
$x(1), \dots, x(g) \in \Gamma({\mathcal{U}}, {\mathcal{A}})$, such that
${\mathcal{A}}|_{{\mathcal{U}}} = {\mathcal{O}}[x(1),\dots, x(g)]$.  By definition of the
${\mathcal{A}}$-module structure on ${\mathcal{F}}$,
\begin{align} &x(i_1) \dots x(i_k)(\sigma \otimes L)\notag\\ = &x(i_1)
\dots x(i_{k-1})(\nabla_{i_k}(\sigma) \otimes L - \sigma\otimes \xi_{i_k}
L)\ .
\end{align} It now follows by induction on $k$ that $\sigma\otimes 1$
generates ${\mathcal{F}}$ freely as a ${\mathcal{A}}$-module.
\end{pf*}

A second proof will be given in section \ref{duality}.

\section{Coherence, Duality and Holonomicity}\label{duality}

\subsection{} If $\G$ is an $\o_Y$-module, we may form
$\d\otimes_{\o} \G$ using right multiplication by functions  as the
$\o_Y$-module structure on $\d$.  This leaves us with left multiplication by
elements of $\d$ to give a left $\d$-module structure.  Similarly, if
$\f$ is an $\o_X$-module, we may form the ${\mathcal{A}}$-module
${\mathcal{A}}\otimes\f$. Note that
\begin{align}\label{adjoint} Hom_{\d}(\d\otimes_{\o} \G,\ \cdot\
)&=Hom_{\o}( \G,\ \for(\cdot)\ )\ ;
\label{adjoint1}\\  Hom_{\b}(\b\otimes\f,\ \cdot\ )&=Hom_{\o}( \f,\
\for(\cdot)\ )\ ,
\label{adjoint2}\end{align} where $\for$ is the forgetful functor.

\begin{Thm}\label{tensor}
\begin{equation}RS_2(\d\overset{L}{\underset{=}{\otimes}}(\cdot))=
\b\overset{L}{\underset{=}{\otimes}}R{\mathcal{S}}_2(\cdot)\ .
\end{equation}
\begin{equation}RS_1(\b\overset{L}{\underset{=}{\otimes}}(\cdot))=
\d\overset{L}{\underset{=}{\otimes}}R{\mathcal{S}}_1(\cdot)\ .
\end{equation}
\end{Thm}
\begin{pf}   By theorem \ref{equiv2}, it suffices to prove the first
equality. Let $\f_1$ be an object in $D^b\text{Mod}(\b)$ and let
$\f_2$ be an object in $D^b\text{Mod}(\o_X)$.  Using \eqref{adjoint1},
\eqref{adjoint2}, Mukai's theorem and theorem \ref{equiv},
\begin{align} Hom(\b&\overset{L}{\underset{=}{\otimes}}\f_2,\f_1)\notag\\
&=Hom(\f_2,\for(\f_1))\notag\\ &=Hom(R\cal
S_1(\f_2),\for(RS_1(\f_1)))\notag\\
&=Hom(\d\overset{L}{\underset{=}{\otimes}}R{\mathcal{S}}_1(\f_2),RS_1(\f_1))\notag\\
&=Hom(RS_2(\d\overset{L}{\underset{=}{\otimes}}R\cal
S_1(\f_2)),T^{-g}(\f_1))\notag\\
&=Hom(RS_2(\d\overset{L}{\underset{=}{\otimes}}T^gR{\mathcal{S}}_1(\f_2)),\f_1)\ .
\end{align} Writing $\f_2$ as $R{\mathcal{S}}_2(\G)$ for some $\G$, we get
\begin{equation} Hom(\b\overset{L}{\underset{=}{\otimes}}R{\mathcal{S}}_2(\G),\f_1)=
Hom(RS_2(\d\overset{L}{\underset{=}{\otimes}}\G),\f_1)\ .
\end{equation} The theorem is proved.
\end{pf}
\noindent{\bf Remark}\  This result appears to conflict with the interchange
of tensor and pontrjagin product stated in theorem Mukai2.   Note, however,
that in the definition of $\d\otimes_{\o} \G$, one is  using
the right $\o$-module structure on $\d$  to define the tensor product and
the left $\o$-module structure on the resulting sheaf.
\vskip 12pt

We digress to give  a
\begin{pf*}{Short proof of Proposition \ref{fm of B}} By theorem
\ref{tensor},
\begin{align} RS_2(\d_{\{0\}\rightarrow Y}) &=RS_2(\d\otimes k(0))\notag\\
&=\b\otimes R{\mathcal{S}}_2(k(0))\notag\\ &=\b\otimes \o_X=\b\ .
\end{align}
\end{pf*}

A more important corollary is the following.  Let
$D^b_{coh}\text{Mod}({\mathcal{A}})$
(resp. $D^b_{coh}\text{Mod}(\d_Y)$) denote the subcategory of complexes with
${\mathcal{A}}$-coherent (resp. $\d$-coherent) cohomology.
\begin{Thm}  The functors  $RS_1$ and $RS_2$ restrict to equivalences
between the categories $D^b_{coh}\text{Mod}({\mathcal{A}})$ and
$D^b_{coh}\text{Mod}(\d_Y)$.
\end{Thm}

\begin{pf}  We give the proof in one direction, the other direction being
the same.   The category $D^b_{coh}\text{Mod}({\mathcal{A}})$ of complexes with
coherent cohomology is  generated by sheaves of the form ${\mathcal{A}}\otimes\f$,
where $\f$ is a coherent $\o_X$-module.  It is well-known that $R\cal
S_1(\f)$  the belongs to $D^b_{coh}\text{Mod}({\mathcal{O}}_Y)$.  By theorem
\ref{tensor},
\begin{equation}
RS_1(\b\otimes\f)=
\d\overset{L}{\underset{=}{\otimes}}R{\mathcal{S}}_1(\f)\ .
\end{equation}
This completes the proof, for
$\d\overset{L}{\underset{=}{\otimes}}(\cdot)$\  maps
$D^b_{coh}\text{Mod}({\mathcal{O}}_Y)$ to $D^b_{coh}\text{Mod}({\mathcal{D}}_Y)$.
\end{pf}

\subsection{} Given a complex $\f \in \text{Ob}\ D^b_{coh}\text{Mod}(\d_Y)$,
its dual complex is
\begin{equation}
\Delta^{\d_Y}(\f) = R {\mathcal{H}}om_{\d_Y}(\f, T^g(\d_Y))\ .
\end{equation} Note that $\Delta^{{\mathcal{D}}_Y}(\f)$ is naturally a complex of
right
$\d$-modules, but we regard it as a complex of left $\d$-modules, using the
antiinvolution
\begin{align}
\d &\longrightarrow \d\notag\\ L = \sum f_I \xi^I &\mapsto \sum(-1)^{|I|}
\xi^I f_I
\ \ .
\end{align}

Given  $\f \in \text{Ob}\ D^b\text{Mod}({\mathcal{A}})$, define
$\Delta^{{\mathcal{A}}}\f = R {\mathcal{H}}om(\f, T^g({\mathcal{A}}))$.
In particular,
\begin{align}
\Delta^{\cal
A}(\b\overset{L}{\underset{=}{\otimes}}\f)
&=\b\overset{L}{\underset{=}{\otimes}}\Delta^{\o_X}(\f)\ ;\\
\Delta^{\d}(\d\overset{L}{\underset{=}{\otimes}}\f)
&=\d\overset{L}{\underset{=}{\otimes}}\Delta^{\o_Y}(\f)\ .
\end{align}

\begin{Prop}\label{grothendieck}
\begin{align} RS_2\Delta^{\d_Y} &= T^{-g}\Delta^{{\mathcal{A}}}RS_2\ ,\\
\Delta^{\d_Y} RS_1 &= RS_1\Delta^{{\mathcal{A}}}T^{-g}\ .
\end{align}
\end{Prop}

\begin{pf}  By theorem \ref{equiv2} it suffices to prove the first equality.
It also suffices to consider the case where $\f$ is an object in
$D^b\text{Mod}(\d_Y)$   of the form $\f = \d \otimes_{\o}\G$, for some
object $\G \in
\text{Ob}\ D^b\text{Mod}(\o_Y)$.  Then
$\Delta^{\d}(\f) = \d \otimes \Delta^{\o_Y}(\G)$. By \cite[p. 161]{Muk}
\begin{equation} R{\mathcal{S}}_2 \Delta^{\o_Y} = T^{-g}\Delta^{\o_X}R{\mathcal{S}}_2\ .
\end{equation} Then
\begin{align} &RS_2\Delta^{\d_Y}(\f) = RS_2(\d \otimes \Delta^{\o_Y}(\G))\\
&= {\mathcal{A}} \otimes R {\mathcal{S}}_2\Delta^{\o_Y}(\G)\ \ \ \ \text{(by theorem
\ref{tensor})}\notag\\  &= {\mathcal{A}} \otimes T^{-g} \Delta^{\o_X} R\cal
S_2(\G)\notag\\ &= {\mathcal{A}} \otimes T^{-g}R {\mathcal{H}}om_{\o}(R{\mathcal{S}}_2(\G),
T^g(\o_X))\notag\\ &= T^{-g}R {\mathcal{H}}om_{{\mathcal{A}}}({\mathcal{A}} \otimes R{\mathcal{S}}_2(\G),
T^g({\mathcal{A}}))\notag\\ &= T^{-g} R {\mathcal{H}}om_{{\mathcal{A}}}(RS_2(\f), T^g({\mathcal{A}})) =
T^{-g}\Delta^{{\mathcal{A}}} RS_2(\f)\ .
\end{align}
\end{pf}

\subsection{}
For the rest of this section we assume $char(k)=0$. Recall that a $\d$-module
$\f$ is said to be holonomic if its  characteristic variety has the least
possible dimension, namely $g$.

\begin{Prop}\label{Borel} Let $\f$ be a
coherent
$\d$-module.  Then $\f$ is holonomic if and only if $H^i(\Delta^{\d}(\f)) =
0$ for $i \ne 0$.
\end{Prop}
\begin{pf} \cite[p. 230]{Bo}\end{pf}

Since holonomicity of a $\d$-module is a local condition,
one expects it to be encoded globally when one takes the Fourier
transform.

\begin{Thm}   Let $\f \in \text{Ob}\ D^b\text{Mod}({\mathcal{A}})$ be
a complex such that the cohomology of $RS_1(\f)$ is concentrated in
a single degree $i$ and  $R^iS_1(\f)$ is $\d$-coherent.
Then $R^iS_1(\f)$ is holomonic if and only if the cohomology of
$RS_1\Delta^{\b}(\f)$ is concentrated in  degree $g-i$.
\end{Thm}

\begin{pf}  Let $\hat{\f}=R^iS_1(\f)$. Regarding $\hat{\f}$ as a
complex in degree $0$,\hfill\break
$H^j(\Delta^{\d}(\hat{\f}))=H^j(\Delta^{\d}T^iRS_1(\f))$.
Then by proposition \ref{Borel}, $\hat{\f}$
is holonomic if and only if
$H^j(\Delta^{\d}T^iRS_1(\f)) = 0$ for $j \ne 0$.
By proposition
\ref{grothendieck},
\begin{align}\label{holonomic}
&H^j(\Delta^{\d}T^iRS_1(\f)) = H^j(\Delta^{\d}RS_1 T^i(\f))\notag\\ &=
H^j(RS_1\Delta^{{\mathcal{A}}}T^{i-g}(\f))= H^j(RS_1T^{g-i}\Delta^{\cal
A}(\f))\ .
\end{align}
This vanishes for $j\ne 0$ if and only if
$R^lS_1\Delta^{{\mathcal{A}}}(\f)$ vanishes for $l\ne g-i$.
\end{pf}
We leave the detailed study of this condition to a future work.

\section{Characteristic Variety}\label{char var}

In many important examples, it is possible to be quite explicit about the
characteristic variety of the transform of a coherent ${\mathcal{A}}$-module.

Let $\{{\mathcal{A}}(m)\}$  be the filtration in example \ref{B itself}.  If $\f$
is a coherent $\b$-module, then one has the usual notion of good filtration
with respect to $\{{\mathcal{A}}(m)\}$. As the \fm\ exchanges local data for global
data, it is worth noting that good filtrations exist globally.


\begin{Prop} Let  $Z$ be a projective scheme of finite type over $k$, and
let ${\mathcal{A}}$ be a sheaf of $\o_Z$-algebras.  Regarding ${\mathcal{A}}$  as an
$\o$-module by letting $\o$ act on the right, assume ${\mathcal{A}}$ is
quasicoherent.  If ${\mathcal{L}}$ is an ample line bundle on $Z$ and ${\mathcal{M}}$ is a
sheaf of coherent left ${\mathcal{A}}$-modules, then there is presentation of
${\mathcal{M}}$ of the form
\begin{equation}\label{presentation} ({\mathcal{A}}\otimes {\mathcal{L}}^{n_2})^{r_2}\too
({\mathcal{A}}\otimes {\mathcal{L}}^{n_1})^{r_1}\too{\mathcal{M}}\too 0\ .
\end{equation}
\end{Prop}
\begin{pf} The proof is the same as for the case ${\mathcal{A}}=\o$.\  (cf.
\cite[p. 122]{H} )
\end{pf}  In particular,
\begin{Cor}  If $\f$ is a sheaf of coherent ${\mathcal{A}}$ (resp. ${\mathcal{D}}_Y$)
modules, then ${\mathcal{M}}$  has a global good filtration by coherent $\o_X$
(resp. $\o_Y$) modules.
\end{Cor}

Let $\f$ be a coherent $\b$-module, and let $\{\f_m\}$ be a good
filtration.  It follows from \eqref{filtration} that
\begin{equation}
Gr\b=Sym({\mathfrak{g}})\otimes\o_X\ .
\end{equation}
Thus, for each $\xi\in{\mathfrak{g}}$ we have a homogeneous map of degree $1$
\begin{equation}
Gr\f\overset{\dot\xi}\too Gr\f\ .
\end{equation}
Identifying $\b$ with $\o[x(1),...,x(g)]$ on a sufficiently small
open set, we get a commutative diagram
\begin{equation}\label{comm}
\begin{CD}
\f_m @>{x(i)}>> \f_{m+1} \\ @VVV @VVV \\ Gr_m\f @>{\dot{\xi}_i}>>
Gr_{m+1}\f
\end{CD} \ \ .
\end{equation} Set $R^jS_1(\f,\psi)=(R^j{\mathcal{S}}_1(\f),\nabla^j)$, and set
$\nabla^j=\sum \omega^i\nabla^j(i)$.  It follows from the explicit formula
for $S_1$ that for all $j$, the diagram
\begin{equation}\label{comm transform}
\begin{CD} R^j{\mathcal{S}}_1(\f_m) @>{\nabla^j(i)} >> R^j{\mathcal{S}}_1(\f_{m+1}) \\
@VVV @VVV \\ R^j{\mathcal{S}}_1(Gr_m\f) @> {R^j{\mathcal{S}}_1(\dot{\xi}_i)}>> R^j\cal
S_1(Gr_{m+1}\f)
\end{CD}
\end{equation} commutes.

Let ${\mathcal{K}}_m$ denote the kernel of the natural map
\begin{equation}
\g\otimes Gr_m\f\too Gr_{m+1}\f\ .
\end{equation}
Let us say that the filtered $\b$-module $\f$ satisfies {\it filtered W.I.T
with index $i$} (cf. \cite[p.156]{Muk}) if
$R^j{\mathcal{S}}_1(\f)=0$ for $j\ne i$ and the same is also true of $\f_m$ and
${\mathcal{K}}_m$ for $m$ sufficiently large. Following Mukai, we denote the
surviving cohomology sheaf  of $RS_1(\f)$ by
$\hat{\f}$.  Then we have
\begin{equation}
0\too \hat{{\mathcal{K}}}_m\too \g\otimes \widehat{Gr_m\f}\too
\widehat{Gr_{m+1}\f}\too 0
\end{equation} for $m$ sufficiently large.
By the preceding remarks, we therefore have
\begin{Prop}\label{good filtration}  Let $\f$ be a filtered
$\b$-module satisfying filtered W.I.T..  Then

1. $\{\hat{\f_m}\}$ is a good $\d$-filtration on $\hat{\f}$.

2. $\widehat{Gr_m\f}=Gr_m\hat{\f}$ for $m$ sufficiently large.
\end{Prop}

Assume $\f$ satisfies filtered W.I.T. with index $i$, and let
${{\mathcal{I}}}(\hat{\f})\subset Sym({\mathfrak{g}})\otimes\o_Y$ denote the characteristic
ideal sheaf of
$\hat{\f}$. Fix an affine open subset ${\mathcal{U}}\subset Y$, let $A=\Gamma(\cal
U,\o_Y)$ and let $I=
\Gamma({\mathcal{U}},{\mathcal{I}}(\hat{\f}))$.  This ideal may be described as follows.
We have a map
\begin{equation}
\text{Sym}({\mathfrak{g}}) \overset{gr\psi}{\longrightarrow} H^0(X,{\mathcal{E}}nd (Gr\cal
F))
\end{equation}
coming from the $Gr\b$-module structure on $Gr\f$.  Redefining the
filtration if necessary, we may assume that ${\mathcal{K}}_m$ and $\f_m$ satisfy
W.I.T. for all $m$.   Since ${\mathcal{U}}$ is affine,
proposition \ref{good filtration} gives
\begin{align}
\Gamma({\mathcal{U}},Gr\hat{\f})&=H^i(X\times\cal
U,\p\otimes\alpha_1^*(Gr\f))\notag\\ &=H^i(X,\alpha_{1*}(\p|_{X\times\cal
U})\otimes(Gr\f))\ .
\end{align}
It is clear from the construction that the
$\text{Sym}(\g)$-module structure on
$\Gamma({\mathcal{U}},Gr\hat{\f})$ is given by the composition
\begin{align}
\text{Sym}({\mathfrak{g}})\overset{gr\psi}{\longrightarrow}
H^0(X,{\mathcal{E}}nd (Gr{\mathcal{F}}))\rightarrow H^0(X,&\cal
End(Gr{\mathcal{F}} \otimes
\alpha_{1*}(\p|_{X\times{\mathcal{U}}})))\notag\\  &\rightarrow
\text{End}(H^i(X, Gr{\mathcal{F}} \otimes \alpha_{1*}(\p|_{X\times{\mathcal{U}}})))\ .
\end{align}  Putting this together with the $A$-module structure,
we have a map
\begin{equation}   A\otimes \text{Sym}({\mathfrak{g}})
\overset{\lambda_{{\mathcal{U}}}}{\longrightarrow}
\text{End}(\Gamma({\mathcal{U}},Gr\hat{\f}))\ .
\end{equation} Thus we have
\begin{Thm}  \label{char}
\begin{equation}
\Gamma({\mathcal{U}},{\mathcal{I}}(\hat{\f}))=\sqrt{\text{ker}(\lambda_{{\mathcal{U}}})}\ .
\end{equation}
\end{Thm}
We will make use of this result in section \ref{nak}.

\section{The Krichever Construction}\label{krich}

Let us briefly explain how the Krichever construction fits into the present
framework.  Assume now that $X=Y=Jac(C)$, where $C$ is a smooth curve of
positive genus.  Pick a base point $P\in C$ and let $C\overset a\too X$ be
the associated abel map.  Taking $\f=a_*(\o_C(*P))$,  it is easy to see
that $\f$ admits a $\b$-module structure.  It suffices to take
representative cocycles $\{c_{nm}(i)\}$ for a basis of $H^1(C,\o)$ with
respect to an open cover $\{{\mathcal{U}}_n\}$, and solve the equations
\begin{equation}
c_{nm}(i)=f_n(i)-f_m(i)
\end{equation}
with $f_n(i)\in\Gamma({\mathcal{U}}_n,\o(*P))$.  (Note that we are identifying
$\o(*P)$ with its endomorphism sheaf.)   An important ingredient here is
the  {\it flag} on $H^1(C,\o)$ coming from the natural maps
\begin{equation}
H^0(C,\o(nP)/\o)\too H^1(C,\o)\ .
\end{equation}
Denoting the images of these maps by $V_n$, we have a sequence of
subsheaves of $\d$,
\begin{equation}
\o_Y=\d_0\subset\d_1\subset ...\ ,
\end{equation}
where $\d_n$ is the $\o$-algebra generated by the vector fields belonging
to $V_n$.  In particular,
\begin{equation}
\d_1=\o[\xi]\ ,
\end{equation}
where $\xi$ is a basis of the one-dimensional space $V_1$.  Now we may also
filter $\b$ by subalgebras
\begin{equation}
\o_X=\b_0\subset\b_1\subset ...\
\end{equation}
in the same way, so that if we set
\begin{equation}
\b_i(m)=\b_i\cap\b(m)\ ,
\end{equation}
then
\begin{equation}
Gr_i\b=Sym(V_i)\otimes\o\ .
\end{equation}
Theorem \ref{equiv2} extends to this more general situation:
$$\text{The \fm\ gives an equivalence of categories}$$
$$D^b\text{Mod}(\b_i) \leftrightarrow D^b\text{Mod}(\d_i)\ .$$
Now there is an essentially canonical $\b_1$-module structure on $\o(*P)$.
Indeed, let $z$ be a local parameter at $P$.  Then $a^*(\b_1)$ is the
subsheaf of $\o(*P)[x]$ characterized by the following growth condition:
If ${\mathcal{U}}$ is a neighborhood of $P$, then
\begin{equation}
\Gamma({\mathcal{U}},a^*(\b_1))=\{ \sum f_ix^i \in
\Gamma({\mathcal{U}},\o(*P)[x])\ |
\ \text{for all}\ j, \sum\limits_{i\ge j} \binom ij \frac{f_i}{z^{i-j}}\in
\Gamma({\mathcal{U}},\o)\ \}\ .\end{equation}
Then $\o(*P)$ is in fact a sheaf of $\b_1|_C$-algebras under
\begin{align}
\b_1|_C&\too \o(*P)\notag\\
f(x)&\mapsto f(x=0)\ .
\end{align}
The key to the Krichever construction is
\begin{Prop}  Let $\G=\widehat{\o(*P)}$ regarded as $\d_1$-module.  Then
$\G|_{X-\Theta}$ is {\it canonically} isomorphic to $\d_1$.
\end{Prop}
This proposition is simply a translation into the language of this paper of
well-known results.   The canonical generator is the Baker-Akhiezer
function.  The discussion given here is similar to that of \cite{R}.  The
important point is that
\begin{equation}
H^0(C,\o(*P))=End_{\b_1}(\o(*P))=End_{\d_1}(\G)\ .
\end{equation}
Thus, if we let $\chi$ denote the canonical generator of $\G|_{X-\Theta}$,
then for all $f\in H^0(C,\o(*P))$, there exists a unique $L_f\in
H^0(X-\Theta,\d_1)$ such that
\begin{equation}
f\chi=L_f\chi\ .
\end{equation}
This is the Burchnall-Chaundy representation of $H^0(C,\o(*P))$, done for
all line bundles at once.  The famous result is that the $L_f$'s satisfy
the KP-hierarchy.

\section{Further Examples}\label{nak}

Nakayashiki has studied generalizations of the Krichever construction in
which the curve and point are replaced by an arbitrary variety together
with an ample hypersurface.  In this section we will illustrate the results
of section  \ref{char var} using a somewhat more general version of his
examples.  In particular we will obtain some refinements of his
results about characteristic varieties.

Let $Z$ be a smooth projective variety, and take $X$ to be its albanese
variety.   Assume for simplicity that the albanese map
\begin{equation}
Z\overset a\too X
\end{equation}
is an imbedding.  Let $\b_Z=a^*(\b)$.    Then an
$\b_Z$-module is the same as an $\b$-module supported on $Z$, so we have a
functor
$$\b_Z\text{-modules} \too \d_Y\text{-modules}\ .$$
Let  $D\subset Z$
be an ample  hypersurface.
As in the Krichever construction, $\o(*D)$ may be given the structure of an
algebra over $\b_Z$, making $\widehat{\o(*D)}$ a $\d$-module.    (Note that
$\o(*D)$ is W.I.T. of index 0.)  Nakayashiki refers to such $\d$-modules as
BA-modules.

In general, when one has a splitting on a sheaf ${\mathcal{R}}$, there is an
induced map
\begin{equation}
\g \longrightarrow H^0(X,{\mathcal{E}}nd({\mathcal{R}})/\o)
\end{equation} by virtue of \eqref{cobdry}.  In the present example, this is
a map
\begin{equation}
\g \overset{\psi'}{\longrightarrow} H^0(Z, \o(*D)/\o)
\end{equation} splitting the natural map
\begin{equation} H^0(Z, \o(*D)/\o) \longrightarrow H^1(Z, \o) = \g\ .
\end{equation}
Conversely, let $\g \overset{T}{\rightarrow} H^0(Z, {\mathcal{O}}(rD)/{\mathcal{O}})$ be
any left inverse of the natural map $H^0(Z, {\mathcal{O}}(rD)/{\mathcal{O}})
\rightarrow H^1(Z, {\mathcal{O}}) = \g$.  Then there is a splitting $\psi$ on
${\mathcal{O}}(*D)$ such that  $\psi' = T$.   We say that the splitting is {\it
associated} to
$T$.  Equivalently,
let $D_r$ denote the scheme $(D, {\mathcal{O}}/{\mathcal{O}}(-rD))$ and let ${\mathcal{M}} = \cal
O(rD)/{\mathcal{O}})$, which we regard as a line bundle on $D_r$.  Then we think
of $\psi$ as being associated to the rational map
\begin{equation}
D_r\too \Bbb P(\g^*)
\end{equation}
induced by the linear system $T(\g)$.

Now let ${\mathcal{O}}(*D)$ be endowed with a ${\mathcal{A}}_Z$-module structure
associated to a fixed map $T$.  From $T$ we get a map
\begin{equation}
\g \otimes {\mathcal{O}}_{D_r} \overset{t}{\longrightarrow} {\mathcal{M}}\ ,
\end{equation} from which we may construct a Koszul complex (where ${\mathcal{O}} =
{\mathcal{O}}_{D_r})$
\begin{equation} 0 \rightarrow \wedge^g \g \otimes {\mathcal{M}}^{-g + 1}
\cdots \overset{\wedge_3 t}{\longrightarrow} \wedge^2 \g
\otimes {\mathcal{M}}^{-1} \overset{\wedge_2 t}{\longrightarrow} \g
\otimes {\mathcal{O}} \overset{t}{\longrightarrow} {\mathcal{M}} \rightarrow 0.
\end{equation}

Now let ${\mathcal{H}}$ be a sheaf of coherent ${\mathcal{O}}_Z$-modules, and set
\begin{equation}
{\mathcal{F}} = {\mathcal{H}} \otimes {\mathcal{O}}(*D)\ ,
\end{equation} regarded as a ${\mathcal{A}}_Z$-module.  We make the simplifying
hypothesis that
the injections
\begin{equation}
\o(krD) \longrightarrow \o((k+1)rD)
\end{equation} induce injections
\begin{equation}
\f_k \longrightarrow \f_{k+1}\ ,
\end{equation}
where $\f_k=
{\mathcal{H}} \otimes {\mathcal{O}}(krD)$.
This is equivalent to
\vskip 14pt
\noindent{\bf Assumption}
\begin{equation}
{\mathcal{T}}\kern -.2em or^1({\mathcal{H}},\o_D)=0\ .
\end{equation}   This gives us a filtration on $\f$, with associated graded
sheaf
\begin{equation}
Gr\f_k={\mathcal{H}}|_{D_r}\otimes{\mathcal{M}}^k\ ,
\end{equation}
upon which $\g$ acts through the map $T$.
We want this to be a good filtration.

\begin{Lem}\label{base locus}  The following are equivalent.
\begin{enumerate}
\item $\g \otimes {\mathcal{H}} \rightarrow {\mathcal{M}} \otimes {\mathcal{H}}
\rightarrow 0$ is exact.
\item The baselocus of $T(\g)$ does not meet the support of
${\mathcal{H}}$.
\item The complex $(\wedge^i\g \otimes {\mathcal{M}}^{1-i}) \otimes
{\mathcal{H}}$ is exact
\item If we set ${\mathcal{R}}_j = \text{ker}(\wedge_j t)$, then
\begin{equation} 0 \rightarrow {\mathcal{R}}_j \otimes {\mathcal{H}} \rightarrow \wedge^j\g
\otimes {\mathcal{M}}^{1-j} \otimes {\mathcal{H}} \rightarrow {\mathcal{R}}_{j-1}\otimes
{\mathcal{H}} \rightarrow 0\notag
\end{equation} is exact for all $j$.
\item The filtration $\{\f_k\}$ is a good ${\mathcal{A}}_Z$-filtration.
\end{enumerate}
\end{Lem}

\begin{pf} The equivalence of  1 and 2 follows from Nakayama's lemma.  If 2
holds, then the Koszul complex is exact at points in the support of
${\mathcal{H}}$, from which 3 follows.  Since 1 is part of 3, the first three
statements are equivalent.  Then 4 follows because $\wedge^j
\g \otimes {\mathcal{M}}^{1-j} \rightarrow {\mathcal{R}}_{j-1} \rightarrow 0$ is exact at
all points in $\text{supp}({\mathcal{H}})$, and the
${\mathcal{R}}_j$'s are all projective.  But 1 is part of 4, so 1 through 4 are
equivalent.  Now the filtration ${\mathcal{F}}_k$ is good if  and only if
\begin{equation}\label{also action}
\g \otimes Gr_k{\mathcal{F}} \rightarrow Gr_{k+1}{\mathcal{F}}\too 0
\end{equation} is exact for large $k$.  But $Gr_k{\mathcal{F}} = {\mathcal{H}}|_{{\mathcal{D}}_r}
\otimes {\mathcal{M}}^k$, so 1 and 5 are equivalent.
\end{pf}

Assume now that $\{ {\mathcal{F}}_k\}$ is a good $\cal
A_Z$-filtration.
\begin{Lem} ${\mathcal{F}}$ satisfies filtered W.I.T. with index $0$.
\end{Lem}
\begin{pf}  It is clear that  $\f_k$ satisfies W.I.T. with
index $0$ if $k$ is sufficiently large.  Let
${\mathcal{K}}_k$ denote the kernel of \eqref{also action}.  Since ${\mathcal{M}}$ is
invertible,
\begin{equation}
{\mathcal{K}}_k = {\mathcal{K}}_1 \otimes {\mathcal{M}}^k\ .
\end{equation} Since ${\mathcal{M}}$ is ample on $D_r$, there exists $k$ such that
\begin{equation} H^i(D_r, {\mathcal{K}}_j \otimes {\mathcal{L}}) = 0
\end{equation} for any line bundle ${\mathcal{L}}$ which is the pullback of a
degree-0 line bundle on $X$, any $j \ge k$ and all $i > 0$.
Thus for large $k$,  ${\mathcal{K}}_k$ also satisfies W.I.T. with index $0$.
\end{pf}

This puts us in
the position to apply theorem \ref{char}.
Let ${\mathcal{U}}\subset Y$  be affine open and let $A=\Gamma(\cal
U,\o_Y)$.  Let
$\epsilon_1$ and $\epsilon_2$ denote the projections on $D_r\times Y$. We
must study  the graded
$A\otimes
\text{Sym}(\g)$-module
\begin{equation} H\overset{\text{def}}=
\bigoplus_{k=0}^{\infty}H^0(D_r,{\mathcal{H}}|_{D_r}\otimes
\epsilon_{1*}(\p|_{D_r\times{\mathcal{U}}})\otimes{\mathcal{M}}^k)\ .
\end{equation} Let
\begin{equation} S=H^0(D_r,\bigoplus_{k=0}^{\infty} {\mathcal{M}}^k)\ .
\end{equation} Then $H$ is a graded
$S$-module.   Moreover, we have the map
\begin{equation}\label{linear map}
\g \overset{T}{\longrightarrow} H^0(Z, \o(rD)/\o)=
H^0(D_r, {\mathcal{M}})\ ,
\end{equation}
which induces
\begin{equation}\label{induced}
Sym(\g) \too Sym( H^0(D_r, {\mathcal{M}}))\ .
\end{equation}
Then  $Sym({\mathfrak{g}})$ acts on
$H$  through the composition of \eqref{induced} with the natural
homomorphism
\begin{equation}
Sym( H^0(D_r, {\mathcal{M}}))\too S\ .
\end{equation}
If we apply Proj to the composite map $Sym(\g)\too S$, we recover the
rational map
\begin{equation} D_r\overset{\Psi}{--\rightarrow}\Bbb P(\g^*)\ .
\end{equation}
associated to the linear map $T$, \eqref{linear map}.  Moreover, applying Proj
to the graded $A\otimes S$-module $H$, we get the sheaf
\begin{equation}
\tilde{H}=(\alpha_1^*({\mathcal{H}})\otimes\p)|_{D_r\times{\mathcal{U}}}
\end{equation}
on $D_r\times {\mathcal{U}}$.  Thus, the consideration of $H$  as an $A\otimes
Sym(\g)$-module may be viewed on the sheaf level  as taking the direct
image of $(\alpha_1^*({\mathcal{H}})\otimes\p)|_{D_r\times{\mathcal{U}}}$ under the
map
$\Psi\times 1$.   (Recall that by lemma
\ref{base locus}, $\Psi$ is defined on the support of ${\mathcal{H}}|_{D_r}$.)

The discussion above gives us the main result of this section.
\begin{Thm}\label{charac var}  Let $\G=\hat{\f}$, where $\f=\cal
H\otimes\o_Z(*D)$, ${\mathcal{T}}or^1({\mathcal{H}},\o_D)=0$, and the ${\mathcal{A}}_Z$-module
structure on $\o(*D)$ is associated to $\Psi:D_r--\to \Bbb P(\g^*)$.  Then
the characteristic variety of
$\G$, viewed as a subvariety of $\Bbb P(\g^*)\times Y$, is the support of
the sheaf
\begin{equation}  Gr\G=(\Psi\times 1)_*
((\alpha_1^*(\cal
H)\otimes\p)|_{D_r\times Y})\ .
\end{equation}
\end{Thm}

In \cite{N2} the case $Z=X$ is considered, and it is proved that the
codimension of the characteristic variety is $dim(X)-dim Supp ({\mathcal{H}})$.
The more general theorem \ref{charac var} yields quite detailed information
in this case, and in particular has Nakayashiki's result as a corollary.
We are now dealing with a smooth, ample hypersurface $D\subset X$.  Then
there is a natural class of $\b$-module structures on $\o(*D)$.   These are
described in \cite{N2} in terms of factors of automorphy, but  may also be
seen somewhat more geometrically.  Let ${\mathfrak{h}}$ be the space of vector
fields on $X$.  We have the {\it Gauss map}
\begin{equation}
D\overset\Psi\too \Bbb P({\mathfrak{h}}^*)\ .
\end{equation}
The normal bundle to $D$ is $\Psi^*(\o(1))$, and thus we have a linear map
\begin{equation}
{\mathfrak{h}} \overset{\lambda}{\longrightarrow} H^0(D, {\mathcal{N}})\ .
\end{equation}
The composition of $\lambda$ with the canonical map $H^0(D, {\mathcal{N}})\to\g$ is
an isomorphism.
Thus one may consider ${\mathcal{A}}$-module structures on $\o(*D)$
associated to the Gauss map.  That is, we can choose the $\b$-module
structure  in such a way that the induced map
\begin{equation}
\g \overset{\psi'}{\longrightarrow} H^0(X, \o(*D)/\o)
\end{equation}
is precisely the composition
\begin{equation}
\g\simeq{\mathfrak{h}}\overset{\lambda}{\longrightarrow} H^0(D,{\mathcal{N}})\subset
H^0(X, \o(*D)/\o)\ .
\end{equation}
We will call such an ${\mathcal{A}}$-module structure {\it canonical}.  In the
language of this paper, Nakayashiki's $\d$-modules are obtained as the \fm\
of ${\mathcal{A}}$-modules of the form ${\mathcal{H}}\otimes\o(*D)$, where $\o(*D)$ has
a canonical ${\mathcal{A}}$-module structure.   However, one may as well consider
a more general class of ${\mathcal{A}}$-module structures, namely those
associated to any $g$-dimensional linear system $V\subset H^0(D,{\mathcal{N}})$
which is  basepoint-free and maps isomorphically onto $\g$.  The example
of the Gauss map shows that this is the generic situation.  As a
corollary of  theorem \ref{charac var}, we have
\begin{Thm}  Let  $D\subset X$ be a
smooth ample hypersurface.  Let
$V\subset H^0(D,{\mathcal{N}})$ be a  $g$-dimensional basepoint-free linear
system mapping isomorphically onto $\g$, and let
$\Psi:D\to \Bbb P(\g^*)$ be the corresponding morphism.  Let ${\mathcal{H}}$ be a
coherent
$\o_X$-module such that
${\mathcal{T}}\kern -.2em or^1({\mathcal{H}},\o_D)=0$\ , and set $\G=\widehat{\cal
H\otimes\o(*D)}$. Give $\G$ the $\d$-module structure induced by a
$\b$-module structure on $\o(*D)$ associated to $V$.   Then the
characterisic
variety of $\G$ is
\begin{equation} ss(\G)=\Psi(Supp({\mathcal{H}}|_D))\times Y
\ .\end{equation}
(In particular, $ codim(ss(\G)=dim(X)-dim Supp ({\mathcal{H}})$\ .)
\end{Thm}
\begin{pf} Because ${\mathcal{N}}$ is ample, the morphism $\Psi$ is
finite.   Therefore, tensoring
${\mathcal{H}}|_D$ with a line bundle has no effect on the support of its direct
image.  The result then follows from theorem \ref{charac var}.\end{pf}

\section{Some remarks on commuting  rings of matrix partial differential
operators}\label{PDOs}

In some sense, the origins of the present subject date to work of  Burchnall
and Chaundy on commuting rings of ordinary differential operators
\cite{BC}.  Such rings are always of dimension one.   As Nakayashiki has
observed in \cite{N2}, the \fm\ allows one to represent the ring
$H^0(X,\o(*D))$, $X$ an abelian variety and $D$ a smooth ample hypersurface,
by matrix valued partial differential operators in $g$ variables, the size
of the matrix being the $g$-fold self-intersection number, $D^g$.   We want
to offer some further observations on this question.

Consider again the data $(Z,D,{\mathcal{H}})$ in section \ref{nak}.  Fix an integer
$r$ and a subspace $V\subset H^0(Z,\o(rD)/\o)$ mapping isomorphically onto
its image under the natural map
\begin{equation} H^0(Z,\o(rD)/\o)\too \g\ .
\end{equation} As in section \ref{krich}, we get a subsheaf ${\mathcal{A}}_1\subset
{\mathcal{A}}$ by imitating the construction of ${\mathcal{A}}$, replacing $\g$ by $V$
throughout. Similarly, we have a subsheaf $\d_1\subset \d_Y$ generated over
$\o_Y$ by the vector fields belonging to the image of $V$.    As in
section
\ref{nak}, we have a rational map
\begin{equation} D_r\overset{\Psi}{--\rightarrow}\Bbb P(V^*)\ .
\end{equation} Assuming now that $supp({\mathcal{H}})$ does not meet the baselocus
of $\Psi$, we can introduce a coherent ${\mathcal{A}}_1$-module structure on
$\f={\mathcal{H}} \otimes
\o(*D)$, filtered by the submodules
${\mathcal{H}} \otimes
\o(krD)$ exactly as before.  We now have
\begin{align} Gr_1{\mathcal{A}} &=Sym(V)\otimes\o_Z\\ D_r&=Spec(\o/\o(-rD))\\
{\mathcal{M}}&=\o(rD)/\o\ \text{(thought of as a line bundle on $D_r$)}\\ Gr\f&=
\bigoplus_{l=0}^{\infty}{\mathcal{H}}|_{D_r}\otimes{\mathcal{M}}^l\ ,
\end{align} where the $Sym(V)$-module structure on $Gr\f$ is defined by
the inclusion
\begin{equation}\label{inclusion} V\too H^0(D_r,{\mathcal{M}})\ .
\end{equation} Set
$
{\mathcal{G}}=\hat{\f}$, regarded as a sheaf of $\d_1$-modules. We will examine
conditions under which ${\mathcal{G}}$ is free in a neighborhood of a given line
bundle ${\mathcal{L}}\in Pic^0(Z)$.   In order to have the equality
\begin{equation}
\widehat{Gr\f}=Gr\hat{\f}\
\end{equation} in a neighborhood of ${\mathcal{L}}$, we make the assumptions
\begin{align} &\text{For all}\ k\ge 0\ , i>0\ ,\ H^i(Z,{\mathcal{L}}\otimes\f_k)=0\
.\\ &\text{There exists}\ j\ge 2\ \text{such that for all}\ i\ne j\ ,
\ H^i(Z,{\mathcal{L}}\otimes\f_{-1})=0\ .
\end{align} For  ${\mathcal{G}}$ to be free over $\d_1$ in a neighborhood of $\cal
L$, it is sufficient that the fiber of $Gr{\mathcal{G}}$ at ${\mathcal{L}}$ be free over
$Sym(V)$. We therefore have
\begin{Prop}\label{generically free}    For ${\mathcal{G}}$ to be free over $\d_1$
in a neighborhood of ${\mathcal{L}}$, it is sufficient that
$$
\bigoplus_{l=0}^{\infty}H^0(D_r,({\mathcal{L}}\otimes{\mathcal{H}})|_{D_r}\otimes{\mathcal{M}}^l)
$$ be freely generated as a module over $Sym(V)$.
\end{Prop} Regarding the hypothesis of this proposition, we have  the
following purely algebraic lemma.  Let $M$ be an arbitrary finitely
generated graded
$Sym(V)$-module,  with $M_j=0$ for $j<0$.   For all $j$ we have a complex
\begin{equation}\label{complex}
\wedge^2(V)\otimes M_{j-1}\overset{\beta_j}\too V\otimes
M_j\overset{\alpha_j}\too M_{j+1}\ ,
\end{equation} where $\alpha_j$ is the action of $V$, and $\beta_j$ is
defined by
\begin{equation} v\wedge w\otimes m\mapsto v\otimes wm-w\otimes vm\ .
\end{equation}

\begin{Lem}\label{free}  $M$ is free over $\text{Sym}(V)$ if and only if the
sequence
\eqref{complex} is exact for all $j$.
\end{Lem}

\begin{pf}  If $M = \overset{\ell}{\underset{i=1}{\oplus}}
\text{Sym}(V)[n_j]$, where
\begin{equation}
\text{Sym}(V)[n_i]_{j} = \text{Sym}^{j + n_i}(V)\ ,
\end{equation} then the exactness of the complexes \eqref{complex} follows
from the well-known fact that the natural complex
\begin{equation}
\label{always exact}
\wedge^2(V) \otimes \text{Sym}^k(V) \rightarrow V \otimes
\text{Sym}^{k+1}(V) \rightarrow \text{Sym}^{k+2}(V)
\end{equation} is always exact.

Conversely, suppose that \eqref{complex} is always exact.  For all
$j$, choose a subspace $U_j \subset M_j$ complementary to the image of
$\alpha_{j-1}$.  What we must show is that for all $j$, the natural map
\begin{equation}
\overset{j}{\underset{i=1}{\oplus}} \text{Sym}^i(V) \otimes U_{j-i}
\overset{\gamma_j}{\longrightarrow} M_j
\end{equation} is injective, for then $M$ is isomorphic to
\begin{equation}
\oplus U_j \otimes_k \text{Sym}(V)[-j]\ .
\end{equation} Given $\ell$, assume $\gamma_j$ is injective for $j < \ell$.
Then $M_{\ell-1} \approx \overset{\ell - 1}{\underset{i=0}{\oplus}}
\text{Sym}^i(V) \otimes U_{\ell - 1 - i}$ and we have a commutative diagram
\begin{equation}
\begin{matrix} V \otimes \left(\operatornamewithlimits{\oplus}\limits^{\ell -
1}_{i=0}
\text{Sym}^i (V)
\otimes U_{\ell-1-i}\right)\\
\qquad \qquad \qquad \downarrow \delta \qquad\qquad \searrow
\alpha_{\ell - 1}\\
\operatornamewithlimits{\oplus}\limits^{\ell}_{i=1} \text{Sym}^i (V)
\otimes U_{\ell - i}
\overset{\gamma_{\ell}}{\longrightarrow}M_{\ell}
\end{matrix}
\end{equation} Since
$M_{\ell - 2} \approx \overset{\ell - 2}{\underset{i=0}{\oplus}}
\text{Sym}^i(V) \otimes U_{\ell - 2 - i}$, the exactness of \eqref{always
exact} implies that
$\text{Ker}(\delta) = \text{Im}(\beta_{\ell - 1})$.  Thus
$\text{ker}(\delta) = \text{Ker}(\alpha_{\ell - 1})$, which says that
$\gamma_{\ell}$ is injective.
\end{pf}

The geometric interpretation  of this lemma is the following. Let $S$ be a
scheme over $k$, $V$ a finite dimensional vector space over
$k$,
${\mathcal{M}}$ a line bundle on $S$, and $V
\overset{T}{\longrightarrow} H^0(S, {\mathcal{M}})$ a linear map.  As in section
\ref{nak}, we associate to $T$ a Koszul complex
\begin{equation}
\dots \rightarrow \wedge^3 V \otimes {\mathcal{M}}^{-2} \overset{\wedge_3
t}{\longrightarrow} \wedge^2 V \otimes {\mathcal{M}}^{-1}
\overset{\wedge_2 t}{\longrightarrow} V \otimes {\mathcal{O}}
\overset{t}{\longrightarrow} {\mathcal{M}} \rightarrow 0\ .
\end{equation} Let ${\mathcal{R}}_i = \text{ker}(\wedge_i t)$.  Let ${\mathcal{H}}$ be a
sheaf of ${\mathcal{O}}_S$-modules. Then we have complexes
\begin{equation}
\label{exact1} 0 \rightarrow {\mathcal{R}}_1 \otimes {\mathcal{H}} \rightarrow V \otimes
{\mathcal{H}}
\rightarrow {\mathcal{M}} \otimes {\mathcal{H}} \rightarrow 0
\end{equation}
\begin{equation}
\label{exact2} 0 \rightarrow {\mathcal{R}}_2 \otimes {\mathcal{H}} \rightarrow \wedge^2 V
\otimes {\mathcal{M}}^{-1} \otimes {\mathcal{H}} \rightarrow {\mathcal{R}}_1 \otimes
{\mathcal{H}} \rightarrow 0.
\end{equation} As in lemma \ref{base locus}, these sequences are exact if
and only if
\begin{equation}
\label{exact0} V \otimes {\mathcal{H}} \rightarrow {\mathcal{M}} \otimes {\mathcal{H}} \rightarrow
0
\end{equation} is exact. We therefore have

\begin{Thm}  Assume \eqref{exact0} is exact. Consider the graded
$\text{Sym}(V)$-module
\begin{equation}  M = \underset{j\ge 0}{\oplus} H^0(S,{\mathcal{H}} \otimes \cal
M^j)\ .
\end{equation}  Then $M$ is free over $\text{Sym}(V)$ if and only if
\begin{align}  0 \rightarrow H^1(S, {\mathcal{R}}_2 \otimes {\mathcal{H}} \otimes \cal
M^j)&
\rightarrow \wedge^2 V \otimes H^1(S, {\mathcal{H}} \otimes {\mathcal{M}}^{j-1})
\ \text{is exact for all}\ j\ge 1\ ,\text{and}\\ H^0(S,{\mathcal{R}}_1\otimes\cal
H)=0\ .\end{align}
\end{Thm}
\begin{pf}  The maps $\alpha_j$ and $\beta_j$ in \eqref{complex} are given in
this case by tensoring \eqref{exact1} and
\eqref{exact2} with ${\mathcal{M}}^i$ and taking cohomology:
\begin{equation}
\begin{matrix} 0 \rightarrow H^0(S, {\mathcal{R}}_1 \otimes {\mathcal{H}}\otimes {\mathcal{M}}^j)
\rightarrow V \otimes M_j \overset{\alpha_j}{\longrightarrow} M_{j+1}\\
\qquad \quad \uparrow  \qquad \nearrow
\beta_j\\
\wedge^2(V) \otimes M_{j-1}
\end{matrix}
\end{equation} Therefore, \eqref{complex} is exact for $j\ge 1$ precisely
when
\begin{equation} H^1(S, {\mathcal{R}}_2 \otimes {\mathcal{H}} \otimes {\mathcal{M}}^j) \rightarrow
\wedge^2V \otimes H^1(S, {\mathcal{H}} \otimes {\mathcal{M}}^{j-1})
\end{equation} is injective.  The second condition appears because we have
defined $M$ as a sum over nonnegative $j$, but have not assumed that
$H^0(S,{\mathcal{H}} \otimes {\mathcal{M}}^{-1})=0$.  The theorem then follows from lemma
\ref{free}.
\end{pf}

Let us put all these ingredients together.  We have  a map
\begin{equation} V\overset{T}{\longrightarrow} H^0(D_r, {\mathcal{M}})\ ,  \ \ \cal
M = {\mathcal{O}}(rD)/{\mathcal{O}}\ ,
\end{equation}
giving us subsheaves ${\mathcal{A}}_1\subset {\mathcal{A}}$, $\d_1\subset \d$, and a
${\mathcal{A}}_1$-module structure on $\o(*D)$. We have a sheaf ${\mathcal{H}}$ such that
${\mathcal{T}}or^1({\mathcal{H}}, {\mathcal{O}}_D) = 0$ and
$\text{supp}({\mathcal{H}})$ does not meet the baselocus of
$\text{Im}(T)$.  We have ${\mathcal{G}}=\widehat{{\mathcal{H}}\otimes\o(*D)}$ and
$\f_k={\mathcal{H}}\otimes\o(krD)$.  The
Koszul complex associated to
$T$ is therefore exact on the support of ${\mathcal{H}}|_{D_r}$, and we can apply
the previous theorem in combination with theorem \ref{generically free}.

\begin{Thm}\label{cohomological}  For ${\mathcal{G}}$ to be free over ${\mathcal{D}}_1$ in
a neighborhood of
${\mathcal{L}}
\in \text{Pic}^0(Z)$, the following conditions are sufficient:
\begin{enumerate}
\item  For all $k\ge 0$,\  $i>0$,\ $H^i(Z,{\mathcal{L}}\otimes\f_k)=0$.
\item There exists $j\ge 2$  such that for all $i\ne j$\ ,
$H^i(Z,{\mathcal{L}}\otimes\f_{-1})=0$.
\item $0 \rightarrow H^1(D_r, {\mathcal{R}}_2
\otimes{\mathcal{L}}\otimes {\mathcal{H}} \otimes {\mathcal{M}}^j)
\rightarrow \wedge^2 V \otimes H^1(D_r, {\mathcal{L}}\otimes{\mathcal{H}} \otimes \cal
M^{j-1})$\ is exact for all $j\ge 1$
\item $H^0(D_r,{\mathcal{R}}_1\otimes\cal
L\otimes
{\mathcal{H}})=0$.\end{enumerate}
\end{Thm}

Nakayashiki's embedding may now be recovered.   As in section \ref{nak}, we
take $Z$ to be $X$ itself, and  we assume $D$ is smooth.   We take an $\cal
A$-module structure on $\o(*D)$ associated to a  $g$-dimensional
basepoint-free linear system $V\subset H^0(D,{\mathcal{N}})$
mapping isomorphically onto $\g$.
We take ${\mathcal{H}}=\o(D)$.  This  affects only the
filtration, not the ${\mathcal{A}}$-module structure.  Thus
\begin{equation}
\f_k=\o((k+1)D)\ .
\end{equation}
\begin{Thm} Set ${\mathcal{G}}=\widehat{\o(*D)}$. Then ${\mathcal{G}}$ is a
free
${\mathcal{D}}$-module in a neighborhood of any ${\mathcal{L}} \ne {\mathcal{O}}$.
\end{Thm}
\begin{pf}  By theorem
\ref{cohomological}, it suffices to verify the following:
\begin{enumerate}
\item For all $i>0$,  $k>0$, $H^i(X,{\mathcal{L}}(kD))=0$.
\item There exists $j\ge 2$ such that for all $i\ne j$,
$H^i(X,{\mathcal{L}})=0$.
\item  For all $j\ge 1$,\  $0\to
H^1(D,{\mathcal{R}}_2\otimes {\mathcal{L}}\otimes {\mathcal{N}}^{j+1})\to
\wedge^2(V)\otimes H^1(D,{\mathcal{L}}\otimes {\mathcal{N}}^j)$\ is exact, where $\cal
R_j$ are the kernel sheaves of the Koszul complex $\wedge^j(V)\otimes\cal
N^{1-j}$.
\item $H^0(D, {\mathcal{R}}_1\otimes{\mathcal{L}}\otimes{\mathcal{N}})=0$.
\end{enumerate}
Items 1 and 2 are well-known.  See, for example, \cite[sec. 13, sec.
16]{Mum}.    It then follows from the exact sequence
\begin{equation} 0 \rightarrow {\mathcal{L}}((k-1)D) \rightarrow {\mathcal{L}}(kD)
\rightarrow
{\mathcal{L}} \otimes {\mathcal{N}}^k \rightarrow 0
\end{equation} that for all $k$ and all $0<i<g-1$,\ $H^i(D, {\mathcal{L}} \otimes
{\mathcal{N}}^k) = 0$.     Since ${\mathcal{R}}_{g-1}={\mathcal{N}}^{1-g}$, it follows by
descending induction that
\begin{equation}\label{descent} H^i(D, {\mathcal{R}}_j \otimes {\mathcal{L}} \otimes \cal
N^k) = 0
\text{ for  $0<i<j$ and all $k$.}
\end{equation}
In particular,  3
holds.   Taking $k=1$ we get a stronger statement, also
by descending induction:
\begin{equation} H^i(D, {\mathcal{R}}_j \otimes {\mathcal{L}} \otimes \cal
N) = 0
\text{ for  $i<j$.}
\end{equation}
Thus 4 holds.
\end{pf}

{}From the standpoint of integrable systems, the relevant feature of a $\cal
D$-module is its endomorphism ring.  As we saw in section
\ref{krich}, if $C$ is a curve embedded in its Jacobian,
$\widehat{{\mathcal{O}}(*P)}$ is a ${\mathcal{D}}_1$-module, where ${\mathcal{D}}_1 =
{\mathcal{O}}[\xi]$.  Extending this ${\mathcal{D}}_1$-module structure to a
${\mathcal{D}}$-module structure, we have
\begin{equation} H^0(C,{\mathcal{O}}(*P)) = \text{End}_{{\mathcal{A}}}({\mathcal{O}}_C(*P)) =
\text{End}_{{\mathcal{D}}}(\widehat{{\mathcal{O}}_C(*P)})\ .
\end{equation} It is the analysis of this endomorphism ring which leads to
the $KP$-hierarchy.  Indeed, having trivialized $\widehat{{\mathcal{O}}_C(*P)}$ as a
${\mathcal{D}}_1$-module, its ${\mathcal{D}}$-endomorphisms are then differential operators
in one variable with $g-1$ parameters, satisfying certain nonlinear
differential equations.

More generally, one may hope to associate  dynamics to the endomorphism ring
of a
${\mathcal{D}}$-module coherent over a proper subalgebra, ${\mathcal{O}}[\xi_1,
\dots,
\xi_n]
\subset {\mathcal{D}}$, $n < g$.  Those modules which are free over the smaller
algebra provide a natural starting point for such an investigation.  Note
that the presence of any nontrivial endomorphisms in such a setting is
already a strong condition on the $\d$-module, but one which can easily be
satisfied by the methods presented here. Such examples will be the object of
study in the sequel.

Finally, by way of advertisement, we mention
\vskip 12pt
\noindent{\bf Example:  The Fano Surface.}
\vskip 12pt If $Z \subset \Bbb C \Bbb P^4$ is a smooth cubic hypersurface,
then its family of lines\hfill\break
\hbox{$S = \{\ell \in \text{Gr}(2, 5)\ |\ \ell
\subset Z\}$} is a smooth surface \cite{CG}.  For generic $s \in S$
corresponding to a line $\ell_s$, the set \hbox{$\{t \in S\ |\ \ell_t
\cap
\ell_s = \{pt\}\}$} is a smooth ample hypersurface $D_s$, the incidence
divisor.  We have isomorphisms
\begin{equation}
\text{Alb}(S) \approx \text{Pic}^0(S) \approx J(Z)\ ,
\end{equation} where $J(Z)$ is the intermediate Jacobian of $Z$.  The
albanese map is identified with the assignment $s \rightsquigarrow D_s$,
which gives an embedding
$S
\subset \text{Pic}^0(S)$.  The dimension of $\text{Pic}^0(S)$ is 5.

Setting ${\mathcal{N}} =$ normal bundle of $D_s$ in $S$ and $T_s(S) =$ tangent space
to $S$ at $s$, we have an isomorphism $T_s(S) \approx H^0(D_s, {\mathcal{N}})$.  The
image of the natural map $H^0(D_s, {\mathcal{N}})
\rightarrow H^1(S, {\mathcal{O}})$, gives a subspace
\begin{equation}
\Bbb C^2 \approx V \subset H^1(S, {\mathcal{O}}) \approx \Bbb C^5.
\end{equation}
Set $\d_1=\o[\xi_1,\xi_2]\subset\d$, where $\xi_1$ and $\xi_2$ are a basis
for $V$.   Then the conditions of theorem
\ref{cohomological} are fulfilled with ${\mathcal{H}} = {\mathcal{O}}(2D)$.  Thus
$\widehat{{\mathcal{O}}(*D)}$ is a ${\mathcal{D}}$-module, locally free as a $\cal
D_1$-module at a generic point.  The rank of $\widehat{{\mathcal{O}}(*D)}$ as a
${\mathcal{D}}_1$-module is the degree of the map $D_s \rightarrow \Bbb P^1$
corresponding to the linear system $H^0(D_s, {\mathcal{N}})$.  This degree is also
5.  Thus we have a representation of $H^0(S, {\mathcal{O}}(*D))$ as
$5
\times 5$ matrix partial differential operators in two variables, with
$3 (=
\dim(H^1(S, {\mathcal{O}})) - \dim(V))$ parameters.

\end{document}